\documentclass{amsart}

\usepackage{amsmath,amsthm,amsfonts,amscd,amssymb}

\newtheorem{thm}{Theorem}[section]
\newtheorem{cor}[thm]{Corollary}
\newtheorem{lem}[thm]{Lemma}
\newtheorem{prop}[thm]{Proposition}

\theoremstyle{definition}

\theoremstyle{remark}
\newtheorem{rem}{Remark}[section]

\def\Z{{\mathcal Z}}
\begin{document}

\title{Small zeros of quadratic forms outside a union of varieties}
\author{Wai Kiu Chan, Lenny Fukshansky, and Glenn R. Henshaw}\thanks{The second author was partially supported by a grant from the Simons Foundation (\#208969 to Lenny Fukshansky) and by the NSA Young Investigator Grant \#1210223.}

\address{Department of Mathematics and Computer Science, Wesleyan University, Middletown, CT 06459}
\email{wkchan@wesleyan.edu}
\address{Department of Mathematics, 850 Columbia Avenue, Claremont McKenna College, Claremont, CA 91711}
\email{lenny@cmc.edu}
\address{Department of Mathematics, California State University at Channel Islands, Camarillo, CA 93012}
\email{glenn.henshaw@csuci.edu}

\begin{abstract}
Let $F$ be a quadratic form in $N \geq 2$ variables defined on a vector space $V \subseteq K^N$ over a global field $K$, and $\Z \subseteq K^N$ be a finite union of varieties defined by families of homogeneous polynomials over $K$.    We show that if $V \setminus \Z$ contains a nontrivial zero of $F$, then there exists a linearly independent collection of small-height zeros of $F$ in $V\setminus \Z$, where the height bound does not depend on the height of $\Z$, only on the degrees of its defining polynomials.  As a corollary of this result, we show that there exists a small-height maximal totally isotropic subspace $W$ of the quadratic space $(V,F)$ such that $W$ is not contained in~$\Z$. Our investigation extends previous results on small zeros of quadratic forms, including Cassels' theorem and its various generalizations. The paper also contains an appendix with two variations of Siegel's lemma. All bounds on height are explicit.
\end{abstract}

\maketitle

\def\A{{\mathcal A}}
\def\AA{{\mathfrak A}}
\def\B{{\mathcal B}}
\def\C{{\mathcal C}}
\def\D{{\mathcal D}}
\def\E{{\mathcal E}}
\def\F{{\mathcal F}}
\def\Ff{{\mathfrak F}}
\def\G{{\mathcal G}}
\def\x{{\mathcal H}}
\def\I{{\mathcal I}}
\def\J{{\mathcal J}}
\def\K{{\mathcal K}}
\def\kk{{\mathfrak K}}
\def\L{{\mathcal L}}
\def\LL{{\mathfrak L}}
\def\M{{\mathcal M}}
\def\O{{\mathcal O}}
\def\W{{\omega}}
\def\CC{{\mathfrak C}}
\def\mm{{\mathfrak m}}
\def\MM{{\mathfrak M}}
\def\OO{{\mathfrak O}}
\def\P{{\mathcal P}}
\def\R{{\mathcal R}}
\def\s{{\mathcal S}}
\def\V{{\mathcal V}}
\def\X{{\mathcal X}}
\def\XX{{\mathfrak X}}
\def\Y{{\mathcal Y}}
\def\Z{{\mathcal Z}}
\def\H{{\mathcal H}}
\def\cee{{\mathbb C}}
\def\pee{{\mathbb P}}
\def\que{{\mathbb Q}}
\def\real{{\mathbb R}}
\def\zed{{\mathbb Z}}
\def\hyp{{\mathbb H}}
\def\aaa{{\mathbb A}}
\def\ff{{\mathbb F}}
\def\kk{{\mathfrak K}}
\def\qbar{{\overline{\mathbb Q}}}
\def\kbar{{\overline{K}}}
\def\ybar{{\overline{Y}}}
\def\kkbar{{\overline{\mathfrak K}}}
\def\ubar{{\overline{U}}}
\def\eps{{\varepsilon}}
\def\ahat{{\hat \alpha}}
\def\bhat{{\hat \beta}}
\def\gt{{\tilde \gamma}}
\def\h{{\tfrac12}}
\def\dd{{\partial}}
\def\bfa{{\boldsymbol a}}
\def\bfb{{\boldsymbol b}}
\def\be{{\boldsymbol e}}
\def\bei{{\boldsymbol e_i}}
\def\bff{{\boldsymbol f}}
\def\bc{{\boldsymbol c}}
\def\bm{{\boldsymbol m}}
\def\bk{{\boldsymbol k}}
\def\bi{{\boldsymbol i}}
\def\bl{{\boldsymbol l}}
\def\bq{{\boldsymbol q}}
\def\bu{{\boldsymbol u}}
\def\bt{{\boldsymbol t}}
\def\bs{{\boldsymbol s}}
\def\bfu{{\boldsymbol u}}
\def\bv{{\boldsymbol v}}
\def\bw{{\boldsymbol w}}
\def\bx{{\boldsymbol x}}
\def\bX{{\boldsymbol X}}
\def\bz{{\boldsymbol z}}
\def\bwy{{\boldsymbol y}}
\def\bY{{\boldsymbol Y}}
\def\bL{{\boldsymbol L}}
\def\ba{{\boldsymbol a}}
\def\bb{{\boldsymbol\beta}}
\def\bet{{\boldsymbol\eta}}
\def\bxi{{\boldsymbol\xi}}
\def\bo{{\boldsymbol 0}}
\def\bs{{\boldsymbol S}}
\def\bol{{\boldkey 1}_L}
\def\ep{\varepsilon}
\def\p{\boldsymbol\varphi}
\def\q{\boldsymbol\psi}
\def\rank{\operatorname{rank}}
\def\aut{\operatorname{Aut}}
\def\lcm{\operatorname{lcm}}
\def\sgn{\operatorname{sgn}}
\def\spn{\operatorname{span}}
\def\md{\operatorname{mod}}
\def\Norm{\operatorname{Norm}}
\def\dim{\operatorname{dim}}
\def\det{\operatorname{det}}
\def\Vol{\operatorname{Vol}}
\def\rk{\operatorname{rk}}
\def\ord{\operatorname{ord}}
\def\ker{\operatorname{ker}}
\def\div{\operatorname{div}}
\def\Gal{\operatorname{Gal}}
\def\GL{\operatorname{GL}}
\def\p{\operatorname{p}}
\def\q{\operatorname{q}}
\def\t{\operatorname{t}}
\def\hs{{\hat \sigma}}
\def\chr{\operatorname{char}}

\section{Introduction and statement of results}

The investigation of small-height zeros of quadratic forms was initiated in the celebrated paper of Cassels \cite{cassels:small}, and later continued by a number of authors (see \cite{cassels_overview} for a detailed overview). Let $K$ be a field of characteristic not equal to 2, and let $V$ be an $L$-dimensional subspace of~$K^N$, $N \geq 2$, $1 \leq L \leq N$. Suppose that $F$ is a quadratic form in $N$ variables over $K$ which has a nontrivial zero in $V$.  A natural problem is to determine an upper bound on the infimum of the heights of all nontrivial zeros of $F$ in $V$, whenever an appropriate height function can be defined.  This has first been done by Cassels in the case $K=\que$ and $V=\que^N$, and later generalized to number fields by Raghavan~\cite{raghavan}, to rational function fields by Prestel~\cite{prestel}, and to algebraic function fields by Pfister~\cite{pfister}. Further results included bounds on heights of collections of linearly independent zeros (e.g. Chalk~\cite{chalk}, Schulze-Pillot~\cite{schulze-pillot}) and full isotropic subspaces of the quadratic space~$(V,F)$ (e.g. Schlickewei~\cite{schlickewei}, Schlickewei-Schmidt~\cite{schmidt:schlickewei}, Vaaler~\cite{vaaler:smallzeros, vaaler:smallzeros2}, Fukshansky~\cite{quad_qbar}). Masser \cite{masser} extended Cassels' small-height zero result to zeros of $F$ avoiding a hyperplane, which is equivalent to showing the existence of small-height solutions of a general quadratic diophantine equation. More recently, Masser's result has been extended over number fields and generalized to a finite union of hyperplanes by Fukshansky~\cite{smallzeros} and Dietmann~\cite{dietmann}. The main goal of the present paper is to obtain a result of a similar type, but in a much more general geometric setting when $K$ is a global field: we establish the existence of a linearly independent collection of small-height zeros of $F$ in $V$ outside of an arbitrary finite union of projective varieties defined over $K$, not containing {\it all} zeros of $F$ on~$V$.  {\em Throughout this paper, unless stated the otherwise, $K$ is a global field whose characteristic is not 2}.  Thus, the global field discussed in this paper is either a number field or a function field of transcendence degree one over a finite field $\ff_q$ of $q$ elements, where $q$ is odd.

Let $J \geq 1$ be an integer. For each $1 \leq i \leq J$, let $\s_i$ be a finite set of homogeneous polynomials in $K[X_1, \ldots, X_N]$ and $Z_K(\s_i)$ be its zero set in $K^N$, that is,
$$Z_K(\s_i) = \{\bx \in K^N: P(\bx) = 0 \mbox{ for all } P \in \s_i\}.$$
For the collection $\bs: = \{\s_1, \ldots, \s_J\}$ of finite sets of homogeneous polynomials, define
\begin{equation}
\label{Z_K}
\Z_{\bs} := \bigcup_{i=1}^J Z_K(\s_i),
\end{equation}
and
\begin{equation}
\label{def_M}
M_{\bs} := \sum_{i=1}^J \max \{\deg P: P \in \s_i\}.
\end{equation}

The main result of this paper is the following theorem and its corollary.  The terminology from the algebraic theory of quadratic forms used in their statements will be reviewed below in Section~\ref{notation}, along with definitions of the appropriate height functions.

\begin{thm} \label{miss_hyper} Let $F$ be a nonzero quadratic form in $N$ variables over $K$, $V$ be an $L$-dimensional subspace of $K^N$, and $m$ be the dimension of a maximal totally isotropic subspace of the quadratic space $(V,F)$.  Suppose that, for a collection $\bs$ of finite sets of homogeneous polynomials in $K[X_1, \ldots, X_N]$,  $F$ has a nontrivial zero in $V\setminus \Z_{\bs}$.  Then there exist $m$ linearly independent zeros $\bx_1,\dots,\bx_m$ of $F$ in $V\setminus \Z_{\bs}$ such that
\begin{equation}
\label{h_order_1}
H(\bx_1) \leq H(\bx_2) \leq \dots \leq H(\bx_m),\ h(\bx_1) \leq h(\bx_2) \leq \dots \leq h(\bx_m),
\end{equation}
and for each $1 \leq n \leq m$,
\begin{equation}
\label{miss_hyper_bnd}
H(\bx_n) \leq h(\bx_n) \ll H(F)^{\frac{9L+11}{2}} \H(V)^{9L+12},
\end{equation}
where the implied constant depends only on $K$, $L$, and $M_{\bs}$.
\end{thm}

As a corollary of Theorem~\ref{miss_hyper}, we can also obtain a statement on the existence of nested sequences of totally isotropic subspaces of $(V,F)$ of bounded height not contained in $\Z_{\bs}$.

\begin{cor} \label{miss_subspace} Let $(V,F)$ be the quadratic space of Theorem \ref{miss_hyper}.  Then for each pair of indices $(n,k)$ with $1 \leq n,k \leq m$, there exists a totally isotropic subspace $W^k_n$ of $(V,F)$ such that $\dim_K W^k_n = k$,
$ W^1_n \subset W^2_n \subset \dots \subset W^m_n$, and $W^k_n \nsubseteq \Z_{\bs}$ for each $1 \leq k \leq m$; also
\begin{equation}
\label{cor_bnd}
\H(W^m_n) \ll H(F)^{10L-m+11} \H(V)^{18L+25}.
\end{equation}
In addition, for each $1 \leq k < m$,
\begin{equation}
\label{cor_bnd_1}
\H(W^k_n) \ll  H(\bx_n) \H(W^m_n) \ll H(F)^{\frac{29L + 33 - 2m}{2}} \H(V)^{27L + 37},
\end{equation}
where $\bx_n$ is as in Theorem \ref{miss_hyper}.   All the implied constants depend only on $K$, $L$, $M_{\bs}$, $N$, and $m$.
\end{cor}

We emphasize here that all the implied constants in Theorem~\ref{miss_hyper} and Corollary~\ref{miss_subspace} are effective, and they will be made explicit in the proofs presented in Section~\ref{proofs}. We say a few words about the origins of these constants in Remark~\ref{const} below.

This paper is organized as follows. In Section~\ref{notation} we set the notation, define the necessary constants and height functions, and review the basic terminology in the algebraic theory of quadratic forms. One of the main tools we need in the proof of Theorem~\ref{miss_hyper} and Corollary~\ref{miss_subspace} is a result on the existence of a small-height maximal totally isotropic subspace of a quadratic space. Such a result was obtained over number fields in~\cite{vaaler:smallzeros} and in~\cite{quad_qbar} over $\qbar$. In Section~\ref{function} we prove the function field analog (Theorem~\ref{smallspace}), which provides a key step to the proof of Theorem~\ref{miss_hyper} in the function field case.   A few technical lemmas are proved in Section~\ref{section_miss_hyper}. We then prove Theorem~\ref{miss_hyper}, along with Corollary~\ref{miss_subspace}, in Section~\ref{proofs}.  Finally, Appendix~\ref{null_basis} contains two variations of Siegel's lemma, another tool used in our main argument: one on small-height basis for a vector space over $K$ outside of a finite union of varieties, generalizing the main result of \cite{null} (Theorem~\ref{general_null}), and an ``orthogonal" version of Siegel's lemma for a bilinear space over a function field (Theorem~\ref{orth_siegel}), which is an analogue of Theorem~2.4 of \cite{witt} over a number field and Theorem~6.1 of \cite{quad_qbar} over $\qbar$.

\section{Preliminaries} \label{notation}

\subsection{Notations}
We start with some notation, following \cite{null} and \cite{witt}.   When $K$ is a number field, we write $d = [K:\que]$ for the global degree of $K$ over $\que$, $\D_K$ for its discriminant, $\omega_K$ for the number of roots of unity in $K$, $r_1$ for its number of real embeddings, and $r_2$ for its number of conjugate pairs of complex embeddings; so $d=r_1+2r_2$.

When $K$ is a function field, we fix a $t \in K$ such that $K$ is a finite separable extension of the rational function field $\ff_q(t)$.  Its global degree is $d = [K:\ff_q(t)]$, and the effective degree of $K$ over $\ff_q(t)$ is
$$\mm(K) = \frac{[K:\ff_q(t)]}{[k_0:\ff_q]},$$
where $k_0$ is the algebraic closure of $\ff_q$ in $K$.

Let $M(K)$ be the set of all places of $K$. For each place $v \in M(K)$, we write $K_v$ for the completion of $K$ at $v$, and let $d_v$ be the local degree of $K$ at $v$, which is $[K_v:\que_v]$ in the number field case, and $[K_v:\ff_q(t)_v]$ in the function field case.

If $K$ is a number field, then for each place $v \in M(K)$ we define the absolute value $|\ |_v$ to be the unique absolute value on $K_v$ that extends either the usual absolute value on $\real$ or $\cee$ if $v | \infty$, or the usual $p$-adic absolute value on $\que_p$ if $v|p$, where $p$ is a prime.

If $K$ is a function field, then all absolute values on $K$ are non-archimedean. For each $v \in M(K)$, let $\OO_v$ be the valuation ring of $v$ in $K_v$ and $\MM_v$ the unique maximal ideal in $\OO_v$. We choose the unique corresponding absolute value $|\ |_v$ such that:
\begin{trivlist}
\item (i) if $1/t \in \MM_v$, then $|t|_v = e$,
\item (ii) if an irreducible polynomial $p(t) \in \MM_v$, then $|p(t)|_v = e^{-\deg(p)}$.
\end{trivlist}

\noindent
In both cases, for each non-zero $a \in K$ the following product formula is satisfied:
\begin{equation}
\label{product_formula}
\prod_{v \in M(K)} |a|^{d_v}_v = 1, \quad \mbox{ for all $a \in K^\times$}.
\end{equation}

If $K$ is a number field and $v | \infty$, then for each positive integer $j$ we define, as in \cite{vaaler:smallzeros},
\[ r_v(j) = \left\{ \begin{array}{ll}
    \pi^{-1/2} \Gamma(j/2+1)^{1/j} & \mbox{if $v | \infty$ is real,} \\
    (2\pi)^{-1/2} \Gamma(j+1)^{1/2j} & \mbox{if  $v | \infty$ is complex,}
\end{array}
\right. \]
and
\begin{equation}
\label{constant:B}
B_K(j) = 2|\D_K|^{1/2d} \prod_{v | \infty} r_v(j)^{d_v/d}.
\end{equation}
It is direct to check that $B_K$ is a non-decreasing function of $j$.

If $K$ is a function field, let $g(K)$ be the genus of $K$, and for any integer $\ell \geq 1$ we define
\begin{equation}
\label{EKL}
\E_K(\ell) = e^{\frac{g(K) \ell}{d}}.
\end{equation}
It is well-known that there exists a unique (up to isomorphism) smooth projective curve $Y$ over $\ff_q$ such that $K$ is the field of rational functions on $Y$. We write $n(K)$ for the number of points of $Y$ over $\ff_q$, and $h_K$ for the number of divisor classes of degree zero on $Y$.  For any integer $j \geq 1$, Let
\begin{equation*}
R_K(j) = \frac{n(K)-1}{2} \left( (j-q+2) h_K \sqrt{n(K)} \right)^{\frac{1}{n(K)-1}} + h_K (n(K)-1) \sqrt{n(K)}.
\end{equation*}

We can now define two quantities $C_K(\ell)$ and $A_K(j)$ which will appear in our various height estimates in the subsequent discussion: for any integer $\ell \geq 1$, let
\begin{equation}
\label{CKL}
C_K(\ell) = \left\{ \begin{array}{ll}
\left( \left( \frac{2}{\pi} \right)^{r_2} |\D_K| \right)^{\frac{\ell}{2d}} & \mbox{if $K$ is a number field,} \\
\exp\left(\frac{(g(K)-1+\mm(K))\ell}{\mm(K)}\right) & \mbox{if $K$ is a function field,}
\end{array}
\right.
\end{equation}
and for any integer $j \geq 1$, let
\begin{equation}
\label{AKLM}
A_K(j) = \left\{ \begin{array}{ll}
\left( j \sqrt{2^{r_1} |\D_K|} \right)^{\frac{1}{d}} & \mbox{if $K$ is a number field with $\omega_K \leq j$,} \\
e^{R_K(j)}  & \mbox{if $K$ is a function field with $q \leq j$,} \\
1 & \mbox{otherwise.}
\end{array}
\right.
\end{equation}
Note that $C_K$ and $A_K$ are non-decreasing functions of $\ell$ and $j$, respectively.

\subsection{Heights} \label{heights}

For each $v \in M(K)$, we define a local height $H_v$ on $K_v^N$ by
$$H_v(\bx) = \max_{1 \leq i \leq N} |x_i|^{d_v}_v,\quad \mbox{ for each $\bx \in K_v^N$}.$$
Also, for each $v \mid \infty$ we define another local height
$$\H_v(\bx) = \left( \sum_{i=1}^N |x_i|_v^2 \right)^{d_v/2} \quad \mbox{ for each $\bx \in K_v^N$}.$$
Then we can define two slightly different global height functions on $K^N$:
\begin{equation}
\label{global_heights}
H(\bx) = \left( \prod_{v \in M(K)} H_v(\bx) \right)^{1/d},\ \ \H(\bx) = \left( \prod_{v \nmid \infty} H_v(\bx) \times \prod_{v | \infty} \H_v(\bx) \right)^{1/d},
\end{equation}
for each $\bx \in K^N$. These height functions are {\it homogeneous}, in the sense that they are defined on the projective space over $K^N$, thanks to the product formula (\ref{product_formula}).  It is easy to see that
\begin{equation}
\label{ht_ineq_sqrt}
H(\bx) \leq \H(\bx) \leq \sqrt{N} H(\bx).
\end{equation}
Notice that in case $K$ is a function field, $M(K)$ contains no archimedean places, and so $H(\bx) = \H(\bx)$ for all $\bx \in K^N$. We also define the {\it inhomogeneous} height
$$h(\bx) = H(1,\bx),$$
which generalizes the Weil height on algebraic numbers.  Clearly, $h(\bx) \geq H(\bx)$ for each $\bx \in K^N$.  We extend the height functions $H$ and $h$ to polynomials by evaluating the height of their coefficient vectors.  

If $X$ is a matrix with $\bx_1, \ldots, \bx_L$ as its columns, then $\H(X)$ is always the height $\H(\bx_1\wedge \dots \wedge \bx_L)$, where $\bx_1 \wedge \dots \wedge \bx_L$ is viewed as a vector in $K^{\binom{N}{L}}$ under the standard embedding. More specifically, the coordinates of this vector are determinants of $L \times L$ minors of $X$, and we define
$$H_v(X) = H_v(\bx_1 \wedge \dots \wedge \bx_L)$$
for each non-archimedean place $v \in M(K)$, as well as
$$\H_v(X) = \H_v(\bx_1 \wedge \dots \wedge \bx_L)$$
for each archimedean place $v \in M(K)$. Then
$$\H(X) = \left( \prod_{v \nmid \infty} H_v(X) \times \prod_{v | \infty} \H_v(X) \right)^{1/d}.$$
Let $V$ be an $L$-dimensional subspace of $K^N$.  Then there exist $N \times L$ matrix $X$ and $(N-L) \times N$ matrix $A$, both are over $K$, such that
$$V = \{ X \bt : \bt \in K^L \} = \{ \bx \in K^N : A \bx = 0 \}.$$
The Brill-Gordan duality principle \cite{gordan:1} (also see Theorem 1 on page 294 of \cite{hodge:pedoe}) implies that $\H(X)=\H(A^t)$, and $\H(V)$ is defined to be this common value. Due to the product formula, this global height does not depend on the choice of the basis. This coincides with the choice of heights in \cite{vaaler:siegel}. When $K$ is a function field, $M(K)$ contains only non-archimedean places, and so we will also write $H(V)$ instead of $\H(V)$. Finally, we may also write $H_v(V)$ or $\H_v(V)$ for the local heights of a basis matrix of $V$ when the choice of this matrix is clear.

Let
\begin{equation}
\label{delta}
\delta = \left\{ \begin{array}{ll}
1 & \mbox{if $K$ is a number field}, \\
0 & \mbox{if $K$ is a function field}.
\end{array}
\right.
\end{equation}

We will also need a few technical lemmas detailing some basic properties of heights. The first one bounds the height of a linear combination of vectors.

\begin{lem} \label{sum_height} For $\xi_1,...,\xi_L \in K$ and $\bx_1,...,\bx_L \in K^N$,
$$H \left( \sum_{i=1}^L \xi_i \bx_i \right) \leq h \left( \sum_{i=1}^L \xi_i \bx_i \right) \leq L^{\delta} h(\bxi) \prod_{i=1}^L h(\bx_i),$$
where $\bxi = (\xi_1,...,\xi_L) \in K^L$, and $\delta$ is as in (\ref{delta}) above.
\end{lem}

\noindent
The second one is an adaptation of Lemma 4.7 of \cite{absolute:siegel} to our choice of height functions, using~\eqref{ht_ineq_sqrt}.

\begin{lem} \label{lem_4.7} Let $V$ be a subspace of $K^N$, $N \geq 2$, and let subspaces $U_1,\dots,U_n \subseteq V$ and vectors $\bx_1,\dots,\bx_m \in V$ be such that
$$V = \spn_K \{ U_1,\dots,U_n,\bx_1,\dots,\bx_m \}.$$
Then
$$\H(V) \leq N^{\delta m/2} \H(U_1) \dots \H(U_n) H(\bx_1) \dots H(\bx_m),$$
where $\delta$ is as in (\ref{delta}) above.
\end{lem}

\noindent
The next one is an adaptation of Lemma 2.3 of \cite{witt} to our choice of height functions, using~\eqref{ht_ineq_sqrt}.

\begin{lem} \label{lem_2.3} Let $X$ be a $J \times N$ matrix over $K$ with row vectors $\bx_1,...,\bx_J$, and let $F$ be a symmetric bilinear form in $N$ variables over $K$ (we also write $F$ for its $N \times N$ coefficient matrix). Then
$$\H(X F) \leq N^{3J \delta/2} H(F)^J \prod_{i=1}^J H(\bx_i),$$
where $\delta$ is as in (\ref{delta}) above.
\end{lem}

\noindent
The next one is Lemma 2.2 of \cite{witt} over any global field.

\begin{lem} \label{intersection} Let $U_1$ and $U_2$ be subspaces of $K^N$. Then
\begin{equation}
\label{inter_eq}
\H(U_1 \cap U_2) \leq \H(U_1) \H(U_2).
\end{equation}
\end{lem}

\begin{rem} \label{strong_intersection} It should be remarked that a stronger version of inequality~\eqref{inter_eq} has been produced in~\cite{schmidt-remark} and~\cite{vaaler:struppeck}, specifically:
\begin{equation}
\label{inter_eq-strong}
\H(U_1+U_2) \H(U_1 \cap U_2) \leq \H(U_1) \H(U_2).
\end{equation}
Unfortunately, for our purposes using the stronger inequality~\eqref{inter_eq-strong} instead of~\eqref{inter_eq} does not seem to have an immediate benefit: while the bounds come out considerably more complicated and hard to read, it is not clear how much better they are, since the quantity $\H(U_1+U_2)$ is usually hard to  estimate non-trivially from below.
\end{rem}

\begin{rem}
An important observation is that due to the normalizing exponent $1/d$ in (\ref{global_heights}) all our heights are {\it absolute}, meaning that they do not depend on the number field or function field of definition.  In particular, we can extend their definitions to the algebraic closure of $K$. Lemmas~\ref{sum_height} - \ref{intersection} also hold verbatim with $K$ replaced by the algebraic closure of $K$.
\end{rem}

\subsection{Quadratic Forms}

Here we introduce some basic language of quadratic forms which are necessary for subsequent discussion.  For an introduction to the subject, the readers are referred to, for instance, Chapter 1 of \cite{scharlau}.  For the sake of more generality, we allow $K$ to be any field of characteristic not 2.  We write
$$F(\bX,\bY) = \sum_{i=1}^N \sum_{j=1}^N f_{ij} X_i Y_j$$
for a symmetric bilinear form in $2N$ variables with coefficients $f_{ij} = f_{ji}$ in $K$, and $F(\bX) = F(\bX,\bX)$ for the associated quadratic form in $N$ variables; we also use $F$ to denote the symmetric $N \times N$ coefficient matrix $(f_{ij})_{1 \leq i,j \leq N}$. In particular, we write $H(F)$ for the height $H$ of the matrix $F$ viewed as a vector in $K^{N^2}$, as in our discussion above. Let $V$ be an $L$-dimensional subspace of $K^N$.  Then $F$ is also defined on $V$, and we write $(V,F)$ for the corresponding quadratic space.

A point $\bx$ in a subspace $U$ of $V$ is called singular if $F(\bx,\bwy) = 0$ for all $\bwy \in U$, and it is called nonsingular otherwise. For each subspace $U$ of $(V,F)$, its radical is the set
$$U^{\perp} := \{ \bx \in U : F(\bx, \bwy) = 0\ \forall\ \bwy \in U \},$$
which is the subspace of all singular points in $U$.   We define $\lambda(U):=\dim_K U^{\perp}$, and will write $\lambda$ to denote $\lambda(V)$. A subspace $U$ of $(V,F)$ is called regular if $\lambda(U)=0$.

A point $\bo \neq \bx \in V$ is called isotropic if $F(\bx)=0$, and anisotropic otherwise. A subspace $U$ of $V$ is called isotropic if it contains an isotropic point, and it is called anisotropic otherwise. A totally isotropic subspace $W$ of $(V,F)$ is an isotropic subspace such that for all $\bx,\bwy \in W$, $F(\bx,\bwy)=0$.  If $(V,F)$ is isotropic, then all maximal totally isotropic subspaces of $(V,F)$ contain $V^{\perp}$ and have the same dimension.

If two subspaces $U_1$ and $U_2$ of $(V,F)$ are orthogonal, we write $U_1 \perp U_2$ for their orthogonal sum. If $U$ is a regular subspace of $(V,F)$, then $V = U \perp \left( \perp_V(U) \right)$ and $U \cap  \left( \perp_V(U) \right) = \{\boldsymbol 0\}$, where
\begin{equation}
\label{perp_V}
\perp_V(U) := \{ \bx \in V : F(\bx,\bwy) = 0\ \forall\ \bwy \in U \}
\end{equation}
is the orthogonal complement of $U$ in $V$. Two vectors $\bx,\bwy \in V$ are called a hyperbolic pair if $F(\bx) = F(\bwy) = 0$ and $F(\bx,\bwy) \neq 0$; the subspace $\hyp(\bx,\bwy) := \spn_K \{\bx,\bwy\}$ that they generate is regular and is called a hyperbolic plane. An orthogonal sum of hyperbolic planes is called a hyperbolic space. Every hyperbolic space is regular. It is well known that there exists an orthogonal Witt decomposition of the quadratic space $(V,F)$ of the form
\begin{equation}
\label{decompose}
V = V^{\perp} \perp \hyp_1 \perp\ \dots \perp \hyp_{\W} \perp U,
\end{equation}
where $\hyp_1, \dots, \hyp_{\W}$ are hyperbolic planes, and $U$ is an anisotropic subspace which is determined uniquely up to isometry. The integer $\W$ is called the Witt index of $(V,F)$.  The rank of $F$ on $V$ is $r:=L-\lambda$.


\section{Quadratic forms over function fields}
\label{function}

In this section we establish the following analogue of Vaaler's result \cite{vaaler:smallzeros} on small-height totally isotropic subspaces of a quadratic space over a function field~$K$.  Not only it fills a void in the literature but also establishes the existence of a small-height zero as a corollary, which will be a key initial step of the proof of our main results in the function field case.

\begin{thm} \label{smallspace} Let $K$ be a function field, and let $F$ be a nonzero quadratic form in $K[X_1, \dots, X_N]$. Let $V \subseteq K^N$ be an $L$-dimensional vector space, $1 \leq L \leq N$. Suppose that the quadratic space $(V,F)$ has a totally isotropic subspace of dimension $\ell \geq 1$, where $\ell$ is greater than the dimension of the radical of $(V,F)$. Then there exists a totally isotropic subspace $\A \subseteq V$ of dimension $\ell$ such that
\begin{equation}
\label{iso_bound}
H(\A) \leq q^{(L - \ell)^2 g(K)/d} H(F)^{(L - \ell)/2} H(V).
\end{equation}
\end{thm}

For any $v \in M(K)$, and let $\XX \subseteq K_v^N$ be an $L$-dimensional subspace, $1 \leq L \leq N$.  Let $\bx_1,\dots,\bx_L$ be a basis for $\XX$ and $X = (\bx_1\ \dots\ \bx_L)$ the corresponding $N \times L$ basis matrix. For a subset $I$ of $\{1,\dots,N\}$ of cardinality $L$, define $_IX$ to be the $L \times L$ submatrix of $X$ consisting of the rows indexed by $I$. Let $J=J_v$ be such a subset so that $|\det(_JX)|_v$ is maximal. We define a matrix $\P_v=\P_v(X)$ as in (4.3) of \cite{vaaler:smallzeros}:
\begin{equation}
\label{P_v}
\P_v = X ( _JX)^{-1} \ _J(1_N),
\end{equation}
where $1_N$ is the $N \times N$ identity matrix. As indicated in \cite{vaaler:smallzeros}, $\P_v$ depends only on $\XX$ and not on the choice of a basis matrix $X$, and it acts as a projection operator from $K_v^N$ onto $\XX$ (by left multiplication) which fixes $\XX$ point-wise. With this notation, we can now state some further useful properties of $\P_v$, all can be proved in the same way as their counterparts in \cite{vaaler:smallzeros}.   The readers are reminded that the heights $H$ and $\H$ are the same over a function field, since there are no archimedean places.

\begin{lem} \label{lemmaPQ} Let $Q_v = \frac{1}{2}(1_N + \P_v)$, then:
\begin{enumerate}
\item[(i)] $H_v(\bx) = \max\{H_v(\P_v \bx), H_v((1_N - \P_v)\bx)\}$ for any $\bx \in K_v^N$,
\item[(ii)] $H_v(\bx) = H_v(Q_v \bx)$ for any $\bx \in K_v^N$,
\item[(iii)] $H_v(Y) = H_v(Q_vY)$ for any $N \times L$  matrix $Y$ over $K_v$ with $1 \leq \textnormal{rank}(Y) = L \leq N$.
\end{enumerate}
\end{lem}

\proof
See Lemma 7 (ii), Lemma 8 (iii), and Lemma 9 (iii) of \cite{vaaler:smallzeros}.
\endproof

\begin{lem}\textnormal{\cite[Lemma 4]{vaaler:smallzeros}} \label{proj_lemma}
Let $X$ be a matrix whose columns form a basis of $\XX$, as above, and let $M$ be an integer such that $L < M \leq N$. If $Y$ is an $N \times (M - L)$ matrix over $K_v$ so that $C = (X\ Y)$ is an $N \times M$ matrix of rank $M$, then
$$H_v(C) = H_v(X)H_v((1_N - P_v)Y).$$
\end{lem}

Next, we will need an adelic version of Minkowski's successive minima theorem over function fields, as established in \cite{thunder_mh}. For each $v \in M(K)$, let $\mu_v$ be the Haar measure on $K_v$ such that $\mu_v(\OO_v) = 1$.   Let $K_{\aaa}$ be the ring of adeles of $K$. We choose a Haar measure $\mu_K$ on $K_{\aaa}$ which is given by
$$\mu_K = q^{1 - g(K)} \prod_v \mu_v.$$
We denote the corresponding product measures on $K_v^N$ and $K_{\aaa}^N$ by $\mu_v^N$ and $\mu_K^N$ respectively.

A measurable subset $\s \subseteq K_{\aaa}^N$ is called a {\em coherent system of $\OO_v$-lattices} (see p. 97 of \cite{weil}) if $\s = \prod_v \s_v$ such that each $\s_v$ is an $\OO_v$-lattice on $K_v^N$ and $\s_v = \OO_v^N$ for all but finitely many $v$.  For an $A \in \GL_N(K_{\aaa})$, we define the successive minima of $\s$ with respect to $A$ as in \cite{thunder_mh}:
$$\lambda_i(\s, A) = \inf_{a \in K_{\aaa}^\times} \left\{ \vert a \vert_{\aaa} : \dim_K \left( \spn_K \left( a\s \cap AK^N \right) \right) \geq i \right\},$$
for each $i = 1, \ldots, N$, where $\vert a \vert_{\aaa} := \prod_v \vert a_v \vert_v$ for each $a = (a_v)\in K_{\aaa}^\times$.

\begin{thm} \label{mink_second}
Let $\s$ be a coherent system of $\OO_v$-lattices in $K_{\aaa}^N$, and let $A \in \GL_N(K_{\aaa})$.  Then
$$\mu_{K}^N(\s) \prod_{i = 1}^N \lambda_i(\s, A) \leq q^N \vert \det(A) \vert_{\aaa}.$$
\end{thm}

\proof
This is essentially Corollary 1 of \cite{thunder_mh}, even though only the case $\s = \prod_v \OO_v^N$ is proved there. For each $v \in M(K)$, let $B_v \in \GL_N(K_v)$ be such that $B_v\s_v = \OO_v^N$.  For all but finitely many $v$,  $\s_v = \OO_v^N$, and hence $\vert \det(B_v) \vert_v = 1$.  Thus $B = (B_v) \in \GL_N(K_{\aaa})$.  It is clear that $\lambda_i(\s, A) = \lambda_i(\prod_v \OO_v^N, BA)$ for all $i$.  On the other hand, $\mu_{K}^N(\prod_v \OO_v^N) = \mu_{K}^N(B\s) = \vert \det(B) \vert_{\aaa} \mu_{K}^N(\s)$. The result then follows immediately from Corollary 1 of \cite{thunder_mh}.
\endproof

Next let $F$ be a nonzero quadratic form in $K[X_1, \dots, X_N]$ and let $V \subseteq K^N$ be an $L$-dimensional vector space, $1 \leq L \leq N$. Suppose that the quadratic space $(V,F)$ has a totally isotropic subspace of dimension $\ell \geq 1$.   Among all the totally isotropic subspaces of dimension $\ell$ in $V$, let $\A$ be one with the smallest height (such a subspace exists by the Northcott finiteness property).  For each $v \in M(K)$, let $\P_v = \P_v(\A)$ be the projection of $K_v^N$ onto $\A_v$, the completion of $\A$ at $v$, as defined in \eqref{P_v} above. We will assume throughout that $\ell$ is strictly larger than the dimension of $V^\perp$, the radical of $(V,F)$.  Then $\perp_V(\A)$, defined as in \eqref{perp_V}, is a proper subspace of $V$: otherwise $\A  \subseteq V^{\perp}$, which is not possible by comparing their dimensions.

\begin{prop} \label{propA'}
Let $\bfb$ be a vector in $V \setminus (\perp_V(\A))$ and let $\B$ be the $(\ell+1)$-dimensional subspace of $V$ spanned by $\A$ and $\bfb$. Then there exists an $\ell$-dimensional totally isotropic subspace $\A' \subseteq \B$ such that
\begin{enumerate}
\item[(i)] $\dim(\A \cap \A') = \ell - 1$,
\item[(ii)] $H(\A)^2 \leq H(\A)H(\A') \leq H(F)H(\B)^2$.
\end{enumerate}
In addition, the following inequality is satisfied:
\begin{enumerate}
\item[(iii)] $1 \leq H(F) \prod_{v \in M(K)} H_v((1_N - \P_v)\bfb)^{2/d}$.
\end{enumerate}
\end{prop}

\proof
Since $\bfb \not \in \perp_V(\A)$, the quadratic space $(\B, F)$ is isometric to the orthogonal sum of a hyperbolic plane and a radical of dimension $\ell - 1$.  Therefore, $\B$ has exactly two $\ell$-dimensional totally isotropic subspaces.  One is $\A$, and we call the other $\A'$.  Note that $\A \cap \A'$ is the radical of $\B$, which implies (i).  For each $v \in M(K)$, let $\P'_v$ be the projection of $K_v^N$ onto $\A'_v$, the completion of $\A'$ at $v$, as defined in \eqref{P_v} above.

Now we select a vector $\bwy \in \B \setminus (\A \cup \A')$. Then $\bwy$ is anisotropic and $\bwy = \bfa + \beta\bfb$ for some $\bfa \in \mathcal A$ and $\beta \in K^\times$.  Following the argument in \cite[Page 679]{vaaler:smallzeros2}, we see that equation (4.7) of \cite{vaaler:smallzeros2} implies
$$0 \neq F(\bwy) = 2\beta F(Q_v(1_N - \P_v')\bwy, (1_N - \P_v)\bfb),$$
and hence by Lemma \ref{lemmaPQ}
\begin{eqnarray*}
1 & = & \prod_v \vert (2\beta)^{-1}F(\bwy) \vert_v^{1/d}\\
 & \leq & H(F) \left \{ \prod_v H_v(Q_v(1_N - \P_v')\bwy)\right\}^{1/d}\left\{ \prod_v H_v((1_N - \P_v)\bfb) \right\}^{1/d}\\
 & \leq & H(F) \left \{ \prod_v H_v((1_N - \P_v')\bwy)\right\}^{1/d}\left\{ \prod_v H_v((1_N - \P_v)\bfb) \right\}^{1/d}.
\end{eqnarray*}
Multiplying both sides of the inequality by $H(\mathcal A)H(\mathcal A')$ and applying Lemma \ref{proj_lemma} we obtain (ii) and (iii).
\endproof

\begin{rem} \label{rem_A} Notice that the proof of part (iii) of Proposition~\ref{propA'} holds when $\A$ is any totally isotropic subspace of $(V,F)$. We will use (iii) in this more general context in the proof of Lemma~\ref{fnct_rad_ht} below.
\end{rem}

\begin{prop} \label{enlarge}
Let $U \subset V$ be an $m$-dimensional subspace, $1 \leq m < L$. For each $v \in M(K)$, let $ \P_v(U)$ be the projection of $K_v^N$ onto $U_v$, the completion of $U$ at $v$, as defined in \eqref{P_v} above. There exist $L - m$ linearly independent vectors $\bfb_1, \ldots, \bfb_{L - m}$ in a (vector space) complement of $U$ in $V$ such that
$$\prod_{i = 1}^{L - m} \prod_v H_v((1_N - \P_v(U))\bfb_i)^{1/d} \leq q^{(L - m)g(K)/d} \left( \frac{H(V)}{H(U)} \right).$$
\end{prop}

\proof
Let $Y$ be an $N \times (L - m)$ matrix whose columns form a basis of a (vector space) complement of $U$ in $V$.  For each $v \in M(K)$, let
$$\s_v = \{\bfu \in K_v^{L - m} : (1_N - \P_v(U))Y\bfu \in \OO_v^N \}.$$
The $N \times (L - m)$ matrix $T_v:= (1_N - \P_v(U))Y$ has rank $L - m$.  For, if $(1_N - \P_v(U))Y\bfu = \bo$, then $Y\bfu \in U$, which implies that $\bfu = \bo$. Moreover,
$$H_v(V) = H_v(U)H_v(T_v)$$
by Lemma \ref{proj_lemma}.

Fix $v \in M(K)$.  By rearranging the coordinates if necessary, we may assume that $\left| \det({}_IT_v) \right|_v = H_v(T_v)$ where $I = \{1, \ldots, L - m \}$.  Then the matrix $W_v: = T_v ({}_IT_v^{-1})$ has the identity matrix $1_{L-m}$ on top, and $\vert \det {}_JW_v \vert_v \leq 1$ for all $J \subseteq \{1, \ldots, N\}$ of size $L - m$.  Particularly, this is true for 
$$J = \{1, \ldots, \ell - 1, \ell + 1, \ldots, L-m, L-m+j\}$$
where $1\leq \ell \leq L-m$ and $1 \leq j \leq N - (L - m)$.   However, for this choice of $J$, $\pm \det {}_JW_v$ is the $(\ell, j)$-entry of $W_v$. Thus, all the entries of $W_v$ are in $\OO_v$.  This shows that ${}_IT_v^{-1}\OO_v^{L - m} = \s_v$, and so $\s_v = \OO_v^{L - m}$ for almost all $v$. Then $\s := \prod_v \s_v$ is a coherent system of $\OO_v$-lattices in $K_{\aaa}^{L-m}$.  Moreover,
$$\mu_v^{L - m}(\s_v) = \left| \det({}_IT_v) \right|_v^{-1} = H_v((1_N - \P_v(U))Y)^{-1}.$$
As a result,
$$\mu_K^{L - m}(\s) = q^{(L-m)(1 - g(K))}\prod_v H_v((1_N - \P_v(U))Y)^{-1} = q^{(L-m)(1 - g(K))}\left( \frac{H(U)}{H(V)} \right)^{d}.$$
By Theorem~\ref{mink_second}, the successive minima $\lambda_1, \ldots, \lambda_{L - m}$ of $\s$ with respect to the identity element $(1_{L-m})_{v \in M(K)} \in \GL_{L-m}(K_{\aaa})$ satisfy
$$\lambda_1\cdots \lambda_{L-m} \leq q^{L - m} \mu_K^{L-m}(\s)^{-1} = q^{(L - m)g(K)} \left( \frac{H(V)}{H(U)} \right)^{d}.$$
Let $\{\bfu_1, \ldots, \bfu_{L-m}\}$ be a set of linear independent vectors associated with the successive minima.  In particular, for $i = 1, \ldots, L - m$, we have
$$H_v((1_N - \P_v(U))Y\bfu_i) \leq \lambda_i.$$
The proposition now follows by setting $\bfb_i = Y\bfu_i$ for $i = 1, \ldots, L - m$.
\endproof

\begin{cor} \label{finalcor}
There exists a nonzero vector $\bfb \in V$ such that
\begin{enumerate}
\item[(i)] the subspace $\B := \spn_K \{ \A,\bfb \}  \subseteq V$ has dimension $\ell + 1$,
\item[(ii)] the vectors $(1_N - \P_v)\bfb$, $v \in M(K)$, satisfy
$$\prod_v H_v((1_N - \P_v) \bfb)^{1/d} \leq q^{(L - \ell)g(K)/d} \left( \frac{H(V)}{H(\A)} \right)^{1/(L - \ell)},$$
\item[(iii)] the subspace $\B$ satisfies
$$H(\B) \leq q^{(L-\ell)g(K)/d} H(\A)^{1 - 1/(L-\ell)} H(V)^{1/(L-\ell)}.$$
\end{enumerate}
\end{cor}

\proof
Let $\bfb_1, \ldots, \bfb_{L-\ell}$ be the vectors obtained in Proposition~\ref{enlarge} with $U=\A$ and $m=\ell$, and arrange them so that
$$\prod_v H_v((1_N - \P_v)\bfb_1) \leq \prod_v H_v((1_N - \P_v)\bfb_2) \leq \ldots \leq \prod_v H_v((1_N - \P_v)\bfb_{L-\ell}).$$
Then set $\bfb = \bfb_1$.  Statements (i) and (ii) are then clear, while statement (iii) follows by a direct application of Lemma~\ref{proj_lemma}.
\endproof

\proof [Proof of Theorem~\ref{smallspace}]
Combine Proposition \ref{propA'} (ii) with Corollary \ref{finalcor} (iii).
\endproof

Using the function field version of Siegel's lemma \cite[Corollary 2]{thunder}, we deduce the following results on zeros of $F$ of small height.

\begin{cor} \label{fnct_fld_smallzeros} Let $\A$ be the subspace of $V$ obtained in Theorem~\ref{smallspace}.  There exists a basis $\bx_1, \ldots, \bx_\ell$ of $\A$ such that
\begin{equation}
\label{ffs1}
\prod_{i=1}^\ell H(\bx_i) \leq q^{(L^2 - \ell + \ell^2)g(K)/d} H(F)^{(L-\ell)/2} H(V).
\end{equation}
In particular, if $\ell$ is the dimension of a maximal totally isotropic subspace of $V$, then $V$ contains an isotropic vector $\bfa$ satisfying
\begin{equation}
\label{ffs2}
H(\bfa) \leq q^{(L^2 - \ell + \ell^2)g(K)/\ell d} H(F)^{(L-\ell)/2\ell} H(V)^{1/\ell}.
\end{equation}
\end{cor}

\proof
Inequality \eqref{ffs1} follows by combining Corollary 2 of \cite{thunder} with our Theorem~\ref{smallspace}. Now \eqref{ffs2} follows from \eqref{ffs1} by taking $\A$ to be a maximal totally isotropic subspace of small height  and letting $\bfa$ be the vector of smallest height among $\bx_1, \ldots, \bx_\W$.
\endproof


Under our assumptions on the quadratic space $(V,F)$,  every maximal totally isotropic subspace of $V$ must properly contain the radical of $(V,F)$.  In other words, $(V,F)$ always contains a nonsingular vector.

\begin{cor} \label{smallzero2} There exists a nonsingular isotropic vector $\bfa \in V$ such that
$$H(\bfa) \leq h(\ba) \leq q^{(2L^2-3L+2)g(K)/d}H(F)^{(L-1)/2} H(V).$$
\end{cor}

\proof
In Corollary \ref{fnct_fld_smallzeros}, one of the $\bx_i$'s must be nonsingular.  Then one of the coordinates of $\bx_i$ must be nonzero, say $x_{ij} \neq 0$ for some $1 \leq j \leq N$. Define $\bfa = \frac{1}{x_{ij}} \bx_i$, then $\bfa$ is again a nonsingular zero of $F$ in $V$, one of which coordinates is equal to 1. Hence
$$h(\bfa) = H(\bfa) = H(\bx_i),$$
by the product formula. Then the corollary follows immediately from Corollary \ref{fnct_fld_smallzeros}, as $1 \leq \ell \leq L-1$.
\endproof

We can also produce a bound on the height of the radical of a quadratic space over a function field, which we will use in our main argument.

\begin{lem} \label{fnct_rad_ht} Suppose that the quadratic space $(V,F)$ has rank $1 \leq r < L$. Then
\begin{equation}
\label{sing}
H(V^{\perp}) \leq q^{rg(K)/d} H(F)^{r/2} H(V).
\end{equation}
\end{lem}

\proof
The argument is identical to the proof of Theorem 2 of \cite{vaaler:smallzeros2}, using our Proposition~\ref{enlarge} instead of Theorem 10 of \cite{vaaler:smallzeros}.
\endproof

\begin{rem} \label{rad_height} The analogue of~\eqref{sing} over number fields was established in Theorem~1.3 of \cite{witt}; specifically,
\begin{equation}
\label{sing_height}
\H(V^{\perp}) \leq  B_K(r)^r H(F)^{r/2} \H(V)
\end{equation}
where the constant $B_K(r)$ is defined by \eqref{constant:B} above.
\end{rem}

\section{Technical Lemmas}
\label{section_miss_hyper}


In this section, we establish a few technical lemmas which will be essential in the proofs of Theorem \ref{miss_hyper} and Corollary \ref{miss_subspace}.  We start with a non-vanishing lemma for polynomials, which is a generalization of Theorem 3.1 of \cite{smallzeros}.  For brevity, we will write $\bX$ for the variable vector $(X_1, \ldots, X_N)$ and $K[\bX]$ for $K[X_1, \ldots, X_N]$.

\begin{lem} \label{nonvanish} Let $N, D \geq 1$ be integers and let $P(\bX) \in K[\bX]$ be a nonzero polynomial with degree less than or equal to $D$. Then there exists $\bz \in K^N$ such that $P(\bz) \neq 0$ and
$$h(\bz) \leq A_K(D),$$
where $A_K(D)$ is defined by \eqref{AKLM}.
\end{lem}

\proof
The conclusion of the lemma follows immediately from Lemma 4.1 of \cite{null} combined with the argument in Section~7 of \cite{null} (in particular, see formulas (44) and (45) of \cite{null}, and the paragraph after \cite[Remark 7.2]{null}).
\endproof

We will also need a technical lemma providing an upper bound on the height of a restriction of a polynomial to a subspace.

\begin{lem} \label{restrict} Let $N, D \geq 1$ be integers and let $P(\bX) \in K[\bX]$ be a polynomial of degree $D$. Let $V \subseteq K^N$ be an $L$-dimensional subspace, $1 \leq L \leq N$, such that $P$ is not identically zero on $V$. Let $\bx_1,\dots,\bx_L$ be a basis for $V$ over $K$, write $A$ for the $N \times L$ basis matrix $(\bx_1,\dots,\bx_L)$, and define
$$P_A(Y_1,\dots,Y_L) = P(Y_1\bx_1 + \dots + Y_L\bx_L) \in K[Y_1,\dots,Y_L],$$
so that $P_A$ is a restriction of $P$ to $V$. Then $P_A$ is a polynomial of degree $D$ in $L$ variables over $K$, and
\begin{equation}
\label{ht_PA}
H(P_A) \leq L^{\delta D} H(P) \prod_{i=1}^L h(\bx_i)^D,
\end{equation}
where $\delta$ is as in \eqref{delta}.
\end{lem}

\proof
Notice that
$$P_A(Y_1,\dots,Y_L) = P \left( \sum_{i=1}^L x_{i1}Y_i,\dots,\sum_{i=1}^L x_{iL}Y_i \right),$$
and so for each $v \mid \infty$,
$$H_v(P_A) \leq L^D H_v(P) \max_{1 \leq i \leq L, 1 \leq j \leq N} | x_{ij} |^{D d_v}_v \leq L^D H_v(P) \prod_{i=1}^L H_v(1,\bx_i)^D,$$
while for each $v \nmid \infty$,
$$H_v(P_A) \leq H_v(P) \max_{1 \leq i \leq L, 1 \leq j \leq N} | x_{ij} |^{D d_v}_v \leq H_v(P) \prod_{i=1}^L H_v(1,\bx_i)^D.$$
Then \eqref{ht_PA} follows by taking a product over all places of $K$ while keeping in mind that function fields have no archimedean places.
\endproof

We will also need a lemma on the existence of a small-height hyperbolic pair in a given hyperbolic plane.

\begin{lem} \label{hyper} Let $F$ be a symmetric bilinear form in $2N$ variables over $K$. Let $\hyp \subseteq K^N$ be a hyperbolic plane with respect to $F$. Then there exists a hyperbolic pair $\bx,\bwy$ for $\hyp$ such that
\begin{equation}
\label{x_bound}
H(\bx) \leq h(\bx) \leq \begin{cases}
2 \sqrt{2}\ B_K(1)^2 \,H(F)^{\frac{1}{2}} \H(\hyp) & \mbox{for number field,} \\
q^{4g(K)/d}\, H(F)^{\frac{1}{2}} \H(\hyp), & \mbox{for function field,}
\end{cases}
\end{equation}
as well as
\begin{equation}
\label{y_bound}
H(\bwy) \leq h(\bwy) \leq \begin{cases}
24 \sqrt{2} N^2\ \left( B_K(1) G_K \right)^2\, H(F)^{\frac{3}{2}} \H(\hyp)^3 & \mbox{for number field,} \\
4q^{4g(K)/d} G_K^2\, H(F)^{\frac{3}{2}} \H(\hyp)^3 & \mbox{for function field,}
\end{cases}
\end{equation}
where the constant $G_K$ is $\E_K(2)^{1-\delta} A_K(2) C_K(2)$.
\end{lem}

\proof
The hyperbolic plane $(\hyp,F)$ is a regular 2-dimensional isotropic subspace of $K^N$.  Therefore Corollary 2 of \cite{vaaler:smallzeros} (when $K$ is a number field) and  Corollary~\ref{fnct_fld_smallzeros} above (when $K$ is a function field) imply the existence of $\bo \neq \bx \in \hyp$ such that $F(\bx)=0$ and the height of $\bx$ is bounded as in \eqref{x_bound}. Now Theorem 1.4 of \cite{null} guarantees the existence of a point $\bz \in \hyp$ such that $F(\bz) \neq 0$ and
\begin{equation}
\label{z_bound}
H(\bz) \leq h(\bz) \leq  2 \E_K(2)^{1-\delta} A_K(2) C_K(2) \H(\hyp).
\end{equation}
Since $F(\bx) = 0$ and $F(\bz) \neq 0$, it must be true that $\bx$ and $\bz$ are linearly independent, and hence span $\hyp$. Therefore we must have $F(\bx,\bz) \neq 0$, since $(\hyp,F)$ is regular. Then define
$$\bwy = F(\bz) \bx - 2F(\bx,\bz) \bz.$$
Clearly, $\hyp = \spn_K \{ \bx,\bwy \}$, and it is easy to check that $F(\bwy) = 0$. Once again, regularity of $(\hyp,F)$ implies that $F(\bx,\bwy) \neq 0$, and so $\bx,\bwy$ is a hyperbolic pair for $\hyp$. Finally, we need to produce an estimate on the height of $\bwy$. In case $K$ is a number field, Lemma 2.3 of \cite{smallzeros} implies that
\begin{equation}
\label{y_height}
H(\bwy) \leq h(\bwy) \leq 3N^2 H(F) h(\bx) h(\bz)^2.
\end{equation}
If $K$ is a function field, the argument in the proof of Lemma 2.3 of \cite{smallzeros} implies that
\begin{equation}
\label{y_height_1}
H(\bwy) \leq h(\bwy) \leq H(F) h(\bx) h(\bz)^2,
\end{equation}
since $K$ has no archimedean absolute values. Combining estimates of \eqref{y_height}, \eqref{y_height_1} with \eqref{x_bound} and \eqref{z_bound} produces \eqref{y_bound}.
\endproof

Our next lemma, which works for any field, establishes a basic divisibility property of a polynomial with respect to any fixed monomial ordering.

\begin{lem} \label{poly_order} Let $K$ be any field, and let $P_1(\bX), P_2(\bX) \in K[\bX]$ be two polynomials in $N \geq 1$ variables over $K$. Fix any monomial ordering. Then there exist polynomials $P'_1(\bX),R(\bX) \in K[\bX]$ such that
\begin{equation}
\label{poly_lemma}
P_1(\bX) = P'_1(\bX)+R(\bX)P_2(\bX),
\end{equation}
and the leading monomial of $P_2(\bX)$, with respect to our chosen monomial ordering, does not divide any monomial of $P'_1(\bX)$.
\end{lem}

\proof
Let us write $\LL(P_2)$ for the leading monomial of $P_2(\bX)$ with respect to the chosen monomial order. If $\LL(P_2)$ does not divide any monomial of $P_1$, then set $P'_1(\bX)=P_1(\bX)$ and $R(\bX)=0$, and \eqref{poly_lemma} follows. Hence we may assume that $\LL(P_2)$ divides at least one monomial of $P_1(\bX)$.  Among all such monomials, let $c_{\ba_0} \bX^{\ba_0} := c_{\ba_0} X_1^{a_{01}} \dots X_N^{a_{0N}}$ be the leading one with respect to our chosen monomial order, where $\ba_0 = (a_{01},\dots,a_{0N}) \in \zed_{\geq 0}$ and $c_{\ba_0} \in K$. Define
$$g_1(\bX) = P_1(\bX) - \frac{c_{\ba_0} \bX^{\ba_0}}{\LL(P_2)} P_2(\bX).$$
Now for each $i \geq 1$, let $\ba_i \in \zed_{\geq 0}$ be such that $c_{\ba_i} \bX^{\ba_i}$ is the leading monomial of $g_i(\bX)$ divisible by $\LL(P_2)$. If such monomial exists, define
$$g_{i+1}(\bX) = g_i(\bX) - \frac{c_{\ba_i} \bX^{\ba_i}}{\LL(P_2)} P_2(\bX).$$
Notice that the set of vectors $\ba_i$ as above forms a discrete subset of $\real^N$, which is decreasing with respect to the $L_1$-norm
$$|\bz|_1 := |z_1| + \dots + |z_N|,$$
and is bounded from below by $\bo$, hence it must be a finite set. This means that the process we described terminates, and so there exists some positive integer $k$ such that no monomial of $g_k(\bx)$ is divisible by $\LL(P_2)$. Therefore
$$P_1(\bX) = g_k(\bX) + \sum_{i=1}^{k-1} \frac{c_{\ba_i} \bX^{\ba_i}}{\LL(P_2)} P_2(\bX),$$
and so \eqref{poly_lemma} holds with
$$P'_1(\bX) = g_k(\bX),\ R(\bX) = \sum_{i=1}^{k-1} \frac{c_{\ba_i} \bX^{\ba_i}}{\LL(P_2)},$$
both having the required properties.
\endproof

For the variable vector $\bX$ and any string of increasing indices $I$ from $\{1, \ldots, N\}$, we write ${}_I\!\bX$  for the vector of variables obtained by removing all $X_i$ from $\bX$ whenever $i \in I$.  The notation ${}_I\!K^{N-\vert I\vert}$ denotes the vector space $K^{N-\vert I\vert}$ with coordinates indexed by the indices not in $I$.  For instance, a typical vector $\bw \in {}_1\!K^{N-1}$ is written as $\bw = (w_2, w_3, \ldots, w_N)$.  The next lemma establishes the existence of zeros of especially small height for polynomials of arbitrary degree away from a hypersurface, provided the polynomial in question is of a particular form.

\begin{lem} \label{quad1} Let $N \geq 3$ be an integer, and let $Q(\bX) \in K[\bX]$ be a polynomial of the form
\begin{equation}
\label{Q_form}
Q(\bX) = X_iX_j (c + Q_1({}_{ij}\bX)) + Q_2({}_{ij}\bX)
\end{equation}
for some indices $1 \leq i < j \leq N$, where $0 \neq c \in K$ and $Q_1$, $Q_2$ are polynomials in the $N-2$ variables ${}_{ij}\!\bX$.  Let $P(\bX) \in K[\bX]$ be a polynomial such that there exists $\bo \neq \bz \in K^N$ with $Q(\bz)=0$ and $P(\bz) \neq 0$. Then there exists such a $\bz$ with
\begin{equation}
\label{PQ_height}
H(\bz) \leq h(\bz) \leq A_K(\deg(PQ))^{1 + \deg(Q)} A_K(2\deg(P))^2 H(Q),
\end{equation}
where $A_K$ is the function defined by \eqref{AKLM} above.
\end{lem}

\proof
There is no harm to assume at the outset that $i = 1$ and $j = 2$.  Let us choose a monomial ordering with respect to which the leading monomial of $Q(\bX)$ contains the product $X_1X_2$. Then Lemma \ref{poly_order} guarantees the existence of polynomials $P'(\bX), R(\bX) \in K[\bX]$ such that $P=P'+RQ$ and $X_iX_j$ does not divide $P'(\bX)$. Since for any $\bz \in K^N$ with $Q(\bz)=0$, $P(\bz)=P'(\bz)$, we can assume from the start that $X_1X_2$ does not divide $P(\bX)$, by replacing $P$ with $P'$ if necessary. Then we can write $P(\bX)$ in the form
$$P(\bX) = X_1^k \, G_1({}_2\!\bX) + G_2({}_1\!\bX),$$
for some positive integer $k$ and $(N-1)$-variable polynomials $G_1({}_2\!\bX)$, $G_2({}_1\!\bX)$, where $G_1({}_2\!\bX)$ is either identically zero or has a monomial not divisible by $X_1$.

First assume that $G_1$ is the zero polynomial.  Then $X_1$ does not divide any monomial of $P$, meaning that $P(\bX) =  G_2({}_1\!\bX)$ is a nonzero polynomial in the $N-1$ variables ${}_1\!\bX$.  For the nonzero polynomial $G_2({}_1\!\bX)\bX_2(c + Q_1({}_{12}\!\bX))$,
Lemma \ref{nonvanish} implies that there exists $\bw \in {}_1\!K^{N-1}$ such that
$$G_2(\bw)\, w_2 \, (c+Q_1({}_2\!\bw)) \neq 0,$$
where ${}_2\!\bw: = (w_3, \ldots, w_N)$, with
\begin{equation}
\label{ht_nonzero}
h(\bw) \leq A_K(\deg(P) + \deg(Q_1)+1) \leq A_K(\deg(PQ)).
\end{equation}
Then, for this choice of $\bw$, let
$$z_1 = - \frac{Q_2({}_2\!\bw)}{w_2 (c + Q_1({}_2\!\bw))},$$
and we form a vector $\bz \in K^N$ by putting $z_1$ and $\bw$ together in the obvious manner, that is $\bz = (z_1, w_2, \ldots, w_N)$.  Notice that $Q(\bz)=0$ and $P(\bz) \neq 0$. Moreover,
\begin{eqnarray}
\label{ht_z_1}
h(\bz) & \leq & H(1,z_1) h(\bw)\nonumber\\
& \leq &  H(Q_2({}_{2}\!\bw), z_2 (c + Q_1({}_{2}\!\bw))\, A_K(\deg(PQ)) \nonumber \\
& \leq & h(\bw)^{\deg(Q)} H(Q)\, A_K(\deg(PQ)) \nonumber\\
& \leq & A_K(\deg(PQ))^{1 + \deg(Q)} H(Q).
\end{eqnarray}

Next suppose that $G_1$ is nonzero. Define $r({}_{12}\!\bX)$ to be the sum of all monomials of $G_1({}_2\!\bX)$ that are not divisible by $X_1$. Since  $G_1$ is nonzero, it must have a monomial not divisible by $X_1$, and so $r$ is also a nonzero polynomial. Similar to the argument above, there must exist $\bu \in {}_{12}\!K^{N-2}$ such that  $r(\bu) \neq 0$, and
\begin{equation}
\label{ht_nonzero_1}
h(\bu) \leq A_K(\deg(r)) \leq A_K(\deg(G_1)) \leq A_K(\deg(P) - k).
\end{equation}
Now define
$$g_1(X_1) := G_1(X_1,u_3,\dots,u_N),$$
and
$$g_2(X_2) := G_2(X_2, u_3, \ldots, u_N).$$
Notice in particular that $g_1(X_1)$ is not identically zero, since its constant term $G_1(0)$ is equal to $r(\bu)$, which is nonzero. Let $f_{\bu}(X_1)$ be the function from $K^\times$ into $K^N$ defined by
$$f_{\bu}(X_1): = \left(X_1, - \frac{Q_2(\bu)}{X_1 (c + Q_1(\bu))}, u_{3},\ldots,u_N \right).$$
Then $Q(f_{\bu}(X_1))$ is the zero polynomial.  For simplicity, let $d_1$ and $d_2$ be the degrees of $g_1$ and $g_2$, respectively.  Consider the polynomial
\begin{eqnarray*}
\overline{P}(X_1) & := & X_1^{d_2}\, P(f_{\bu}(X_1)) \\
& = & X_1^{k+d_2} g_1(X_1) + X_1^{d_2} g_2\left( - \frac{Q_2(\bu)}{X_1 (c + Q_1(\bu))} \right).
\end{eqnarray*}
It is direct to check that for any nonzero $x_1 \in K$ $, P(f_{\bu}(x_1)) \neq 0$  whenever $\overline{P}(x_1) \neq 0$.  Notice that $X_1^{k+d_{2}} g_1(X_1)$ is a nonzero polynomial of degree
$$k+d_{1}+d_{2} > d_{2},$$
since $k > 0$, and $X_1^{d_{2}} g_2\left( - \frac{Q_2(\bu)}{X_1 (c + Q_1(\bu))} \right)$ is a polynomial of degree $d_{2}$. Therefore $\overline{P}(X_1)$ is not identically zero, and hence Lemma \ref{nonvanish} implies that there exists $0 \neq \alpha \in K$ such that $\overline{P}(\alpha) \neq 0$ and
\begin{equation}
\label{ht_alpha}
h(\alpha) \leq A_K(\deg(\overline{P})) \leq A_K(k+d_{1}+d_{2}) \leq A_K(2 \deg(P)).
\end{equation}
Now take $\bz = f_{\bu}(\alpha)$, then we have $Q(\bz) = 0$, $P(\bz) \neq 0$, and
\begin{eqnarray}
\label{ht_z_2}
h(\bz) & = & H \left( 1, \bu, \alpha, \frac{Q_2(\bu)}{\alpha (c + Q_1(\bu))} \right) \nonumber \\
& \leq & h(\bu) h(\alpha)\, H \left( 1, \frac{Q_2(\bu)}{\alpha (c + Q_1(\bu))} \right) \nonumber \\
& = & h(\bu) h(\alpha) \,H \left( Q_2(\bu), \alpha (c + Q_1(\bu)) \right) \nonumber \\
& \leq & h(\bu)^{1 + \deg(Q)} h(\alpha)^2\, H(Q) \nonumber \\
& \leq & A_K(\deg(P) - k)^{\deg(Q)+1} A_K(2 \deg(P))^2 \, H(Q),
\end{eqnarray}
where the last inequality follows by combining \eqref{ht_nonzero_1} and \eqref{ht_alpha}.

Now \eqref{PQ_height} follows by combining \eqref{ht_z_1} and \eqref{ht_z_2}, and this finishes the proof of the lemma.
\endproof

\section{Proof of Main Results} \label{proofs}

For any positive integers $\ell$ and $j$, define a number $T_K(\ell, j)$ by
\begin{equation}
\label{T_const_nf}
T_K(\ell,j) = 3^3 2^{\frac{21\ell-21}{2}} \ell^{\frac{27\ell+51}{2}} C_K(\ell)^{9\ell+14} B_K(\ell-1)^{\max\{\ell, 9\}} A_K(j+2)^3 A_K(2j)^2 |\D_K|^{\frac{9}{2d}},
\end{equation}
when $K$ is a number field, and
\begin{equation}
\label{T_const_fff}
T_K(\ell,j) = q^{\frac{(18\ell^2-27\ell+18)g(K)}{d}} C_K(\ell)^{9\ell+15} \E_K(\ell)^{9\ell+15} A_K(j+2)^3 A_K(2j)^2,
\end{equation}
when $K$ is a function field.  Recall that $d$ is the global degree, $\D_K$ is the discriminant, and $g(K)$ is the genus.  The numbers $\E_K(\ell)$, $C_K(\ell)$, and $A_K(j)$ are defined by (\ref{EKL}), (\ref{CKL}), and (\ref{AKLM}), respectively.  It will be shown later that $T_K(L,M+1)$ will serve as the implied constant in \eqref{miss_hyper_bnd}.

\begin{rem} \label{const} The quantities appearing in the definition of this main constant come from lattice co-volume and point counting estimates used in the proofs of different versions of Siegel's lemma (see~\cite{vaaler:siegel}, \cite{thunder}, \cite{absolute:siegel}, and~\cite{null}). In particular, there are natural analogies between the constants in the number field and function field cases. For instance, the quantities $|\D_K|$ in the number field case and $(g(K)-1+\mm(K))\ell/\mm(K)$ in the function field case, appearing in the definition of $C_K(\ell)$ in~\eqref{CKL}, play the role of co-volume in the adelic versions of Minkowski's convex body theorem over these fields, respectively (see~\cite{thunder} for a discussion of this parallel).
\end{rem}

\begin{prop} \label{miss_hyper_one} Let $F$ be a nonzero quadratic form in $N$ variables over $K$, and $V \subseteq K^N$ be an $L$-dimensional subspace, $1 \leq L \leq N$.  Let $P(\bX) \in K[\bX]$ be a polynomial of degree $D$, and assume that there exists a nontrivial zero $\bz$ of $F$ in $V$ such that $P(\bz) \neq 0$. Then there exists such a zero $\bz$ of $F$ with
\begin{equation}
\label{z_bound_miss}
H(\bz) \leq h(\bz) \leq  T_K(L,D)\, H(F)^{\frac{9L+11}{2}}\, \H(V)^{9L+12}.
\end{equation}
\end{prop}

\proof
Let $r$ be the rank of $(V,F)$, and $\lambda$ be the dimension of the radical $V^\perp$; so $r + \lambda$ is $L$, the dimension of $V$ over $K$.
We first handle the case when $K$ is a number field.  First suppose that $P$ is not identically zero on $V^{\perp}$. Then Theorem 1.4 of \cite{null} implies that there exists $\bo \neq \bz \in V^{\perp}$ such that $P(\bz) \neq 0$ and
\begin{equation} \label{hz}
H(\bz) \leq h(\bz) \leq  \lambda A_K(D) C_K(\lambda) \H(V^{\perp}).
\end{equation}
Combining this observation with the upper bound on $\H(V^\perp)$ in \eqref{sing_height}, we obtain
\begin{equation}
\label{mho_1}
H(\bz) \leq h(\bz) \leq
B_K(r)^r \lambda A_K(D) C_K(\lambda) H(F)^{r/2} \H(V).
\end{equation}
Since $F(\bz) = 0$ and the above upper bound for $H(\bz)$ is smaller than the one in (\ref{z_bound_miss}), we are done.

Next assume that $P$ is identically zero on $V^{\perp}$; so $P$ does not vanish at some nonsingular zero of $F$ on $V$.  Let $\bx_1,\dots,\bx_L$ be the small-height basis for $V$, guaranteed by Siegel's lemma (see \cite{vaaler:siegel} and \cite{absolute:siegel} for original results, and Theorems~1.1, 1.2 of \cite{null} for a convenient formulation):
\begin{equation}
\label{siegel_for_V}
\prod_{i=1}^L h(\bx_i) \leq C_K(L) \H(V).
\end{equation}
Let $A = (\bx_1 \dots \bx_L)$ be the corresponding basis matrix and let $F_A$ and $P_A$ be the corresponding restrictions of $F$ and $P$ to $V$, respectively, as defined in Lemma \ref{restrict}. Combining \eqref{ht_PA} with \eqref{siegel_for_V}, we obtain
\begin{equation}
\label{ht_FA}
H(F_A) \leq L^2 C_K(L)^2 \, H(F) \H(V)^2.
\end{equation}
Now notice that for each $\bz \in K^L$, $F_A(\bz) = 0$ and $P_A(\bz) = 0$ if and only if $F(A\bz) = 0$ and $P(A\bz) = 0$, respectively. Moreover, $\bz \in K^L$ is a nonsingular zero of $F_A$ if and only if $A\bz \in V$ is a nonsingular zero of $F$. Also, by Lemma \ref{sum_height} and \eqref{siegel_for_V},
\begin{equation}
\label{ht_A_z}
h(A\bz) = h \left( \sum_{i=1}^L z_i \bx_i \right) \leq L h(\bz) \prod_{i=1}^L h(\bx_i) \leq L C_K(L) h(\bz) \H(V).
\end{equation}
Since $P$ does not vanish at some nonsingular zero of $F$ on $V$, it must be that $P_A$ does not vanish at some  nonsingular zero of $F$ on $K^L$; in particular, the quadratic space $(K^L,F_A)$ must contain a hyperbolic plane. Our next task is to find a hyperbolic pair of bounded height in~$(K^L,F_A)$.

Corollary~1.2 of \cite{smallzeros} guarantees the existence of a nonsingular zero of $\bx \in K^L$ of $F_A$ with
\begin{eqnarray}
\label{ns_x_nf}
h(\bx) & \leq &  2^{\frac{3(L-1)}{2}} L^{\frac{L-1}{2}} |\D_K|^{\frac{1}{2d}} B_K(L-1) H(F_A)^{\frac{L-1}{2}} \nonumber\\
& \leq &  |\D_K|^{\frac{1}{2d}} B_K(L-1) \left( (2L)^{\frac{3}{2}} C_K(L) H(F)^{\frac{1}{2}} \H(V) \right)^{L-1},
\end{eqnarray}
where the last inequality follows by \eqref{ht_FA}.

Let $\bx$ be a nonsingular zero of $F_A$ satisfying \eqref{ns_x_nf}. Then the linear form $F_A(\bx,\bY)$ is not identically zero on $K^L$, and so there must exist a standard basis vector in $K^L$, call it $\bu$, such that $F_A(\bx,\bu) \neq 0$ and $h(\bu) = 1$. Then $\hyp_{xu} := \spn_K \{ \bx,\bu\}$ is a hyperbolic plane in $(K^L,F_A)$ with
\begin{eqnarray}
\label{hyp_xu_nf}
\H(\hyp_{xu})& \leq & L\, H(\bx) H(\bu) \\
 & \leq & (2L)^{\frac{3L-3}{2}} L |\D_K|^{\frac{1}{2d}} B_K(L-1) C_K(L)^{L-1}\, H(F)^{\frac{L-1}{2}} \H(V)^{L-1}, \nonumber
\end{eqnarray}
where the first inequality is given by Lemma~\ref{lem_4.7} and the second follows by \eqref{ns_x_nf}. Let
\begin{equation}
\label{y}
\bwy = F_A(\bu) \bx - 2F_A(\bx,\bu) \bu.
\end{equation}
Then $F_A(\bwy) = 0$ and $F_A(\bx,\bwy) \neq 0$, so $\bx,\bwy$ is a hyperbolic pair. Moreover, \eqref{y_height} and \eqref{y_height_1} state that
$$h(\bwy) \leq  3L^2  H(F_A) h(\bx) h(\bu)^2 = 3L^2 \, H(F_A) h(\bx).$$
Combining this observation with \eqref{ns_x_nf} and \eqref{ht_FA}, we obtain
\begin{equation}
\label{ns_y_nf}
h(\bwy) \leq 3 \times 2^{\frac{3(L-1)}{2}} L^{\frac{3L+5}{2}} |\D_K|^{\frac{1}{2d}} B_K(L-1) \left( C_K(L) H(F)^{\frac{1}{2}} \H(V) \right)^{L+1}.
\end{equation}

Let $U$ be the orthogonal complement of $\hyp_{xu}$ in $(K^L,F_A)$.  It is an $(L-2)$-dimensional space and
\begin{equation*}
U = \left\{ \bv \in K^L : (\bx\ \bu)^t F_A \bv = 0 \right\}.
\end{equation*}
Here we also write $F_A$ for the coefficient matrix of the quadratic form $F_A$. By the Brill-Gordan duality principle discussed in Section~\ref{notation} above, $\H(U)$ is precisely the height $\H$ of the matrix $(\bx\ \bu)^t F_A$, and hence Lemma \ref{lem_2.3} implies that
$$\H(U) \leq L^{3}\, H(F_A)^2 H(\bx) H(\bu),$$
and then \eqref{ht_FA} combined with \eqref{hyp_xu_nf} imply that
\begin{equation}
\label{ht_H_xu_c_1}
\H(U) \leq  2^{\frac{3L-3}{2}} L^{\frac{3L+11}{2}} |\D_K|^{\frac{1}{2d}} B_K(L-1) \left( C_K(L) H(F)^{\frac{1}{2}} \H(V) \right)^{L+3}.
\end{equation}
Let $\bv_1,\dots,\bv_{L-2}$ be the small-height basis for $U$, guaranteed by Siegel's lemma:
\begin{equation}
\label{siegel_for_hyp}
\prod_{i=1}^{L-2} h(\bv_i) \leq C_K(L-2) \, \H(U) \leq C_K(L) \, \H(U).
\end{equation}
Combining \eqref{siegel_for_hyp} with \eqref{ht_H_xu_c_1}, we see that
\begin{equation}
\label{ht_H_xu_c_basis_1}
\prod_{i=1}^{L-2} h(\bv_i) \leq
2^{\frac{3L-3}{2}} L^{\frac{3L+11}{2}} |\D_K|^{\frac{1}{2d}} B_K(L-1) C_K(L)^{L+4} \left( H(F)^{\frac{1}{2}} \H(V) \right)^{L+3}.
\end{equation}
Now define the matrix $B = \left( \bx\ \bwy\ \bv_1 \dots \bv_{L-2} \right) \in \GL_L(K)$, and let
$$Q(\bY) = F_A(B\bY),\ G(\bY) = P_A(B\bY).$$
Then it is easy to see that $Q$ is of the form \eqref{Q_form} (with $i = 1$ and $j = 2$), and so $Q$ and $G$ satisfy the conditions of Lemma~\ref{quad1}. Hence Lemma~\ref{quad1} guarantees the existence of a point $\bw \in K^L$ such that $Q(\bw) = 0$, $G(\bw) \neq 0$, and
\begin{equation}
\label{w_point_ht}
h(\bw) \leq A_K(D + 2)^3 A_K(2 D)^2 H(Q).
\end{equation}
Now notice that standard height inequalities along with \eqref{ht_FA} imply that
\begin{eqnarray*}
H(Q) & \leq & H(B^t F_A B) \leq L^{2} \, H(B)^2 H(F_A) \leq L^{2}\, H(F_A) h(\bx)^2 h(\bwy)^2 \prod_{i=1}^{L-2} h(\bv_i)^2  \\
& \leq & L^{4} C_K(L)^2 \, H(F) \H(V)^2 h(\bx)^2 h(\bwy)^2 \prod_{i=1}^{L-2} h(\bv_i)^2,
\end{eqnarray*}
and so by combining \eqref{ht_H_xu_c_basis_1}  with \eqref{ns_x_nf}, and \eqref{ns_y_nf}, we see that $H(Q)$ is bounded above by
\begin{equation}
\label{G_ht_nf}
3^2 2^{6L-6} L^{9L+17} |\D_K|^{\frac{3}{d}} B_K(L-1)^6 C_K(L)^{6L+9} \left( H(F) \H(V)^2 \right)^{3L+4}.
\end{equation}
Now define $\bz = A (B \bw) \in V$, and notice that $F(\bz) = F_A(B\bw) = Q(\bz) = 0$, and $P(\bz) = P_A(B\bw) = G(\bw) \neq 0$. Hence $\bz$ is precisely the point we are looking for, and to estimate its height first notice that by the same kind of reasoning as in \eqref{ht_A_z},
\begin{equation}
\label{ht_B_w_1}
h(B\bw) = h \left( w_1 \bx + w_2 \bwy + \sum_{i=1}^{L-2} w_{i+2} \bv_i \right) \leq L\, h(\bw) h(\bx) h(\bwy) \prod_{i=1}^L h(\bv_i).
\end{equation}
Then combining \eqref{ht_B_w_1} with \eqref{ht_A_z}, \eqref{w_point_ht}, \eqref{G_ht_nf}, \eqref{ns_x_nf}, \eqref{ns_y_nf}, and \eqref{ht_H_xu_c_basis_1} we obtain \eqref{z_bound_miss}. This completes the proof of the proposition when $K$ is a number field.

Now, let us suppose that $K$ is a function field.  The strategy of the proof for this case is identical to that for the number field case, only at each height estimation we need to replace the height bounds by those specific to the function field case.

When $P$ is not identically zero on $V^\perp$, Theorem 1.4 of \cite{null} still applies.  Combining (\ref{hz}) with Lemma \ref{fnct_rad_ht}, we obtain a nonzero $\bz \in V^\perp$ such that $P(\bz) \neq 0$ and
\begin{equation*}
H(\bz) \leq  h(\bz) \leq q^{\frac{rg(K)}{d}} \lambda A_K(D) C_K(\lambda)\, H(F)^{\frac{r}{2}} \H(V),
\end{equation*}
which is less than the upper bound in (\ref{z_bound_miss}).  This proves the proposition in this special case.

Now, we further assume that $P$ is identically zero on $V^\perp$.  The same argument in the number field case shows that there is a small-height basis $\bx_1, \ldots, \bx_L$ for $V$ such that
$$\prod_{i=1}^L h(\bx_i) \leq C_K(L) \E_K(L)\, \H(V).$$
Let $F_A$ and $P_A$ be the projections of $F$ and $P$, respectively, on $V$ as defined in the number field case, where $A = (\bx_1 \cdots \bx_L)$.  Then Lemma \ref{restrict} implies that
\begin{equation*}
H(F_A) \leq (C_K(L) \E_K(L))^2\, H(F) \H(V)^2,
\end{equation*}
and by Lemma \ref{sum_height},
\begin{equation}\label{fun0}
h(A\bz) \leq C_K(L) \E_K(L) h(\bz)\, H(V)
\end{equation}
for each $\bz \in K^L$.

Corollary \ref{smallzero2} implies that there is a nonsingular zero $\bx \in K^L$ of $F_A$ with
\begin{eqnarray}\label{fun1}
h(\bx) & \leq & q^{\frac{(2L^2 - 3L + 2)g(K)}{d}}\, H(F_A)\nonumber\\
    & \leq & q^{\frac{(2L^2 - 3L + 2)g(K)}{d}} (C_K(L)\E_K(L) H(F)^{\frac{1}{2}} \H(V))^{L-1}.
\end{eqnarray}
Then there exists a standard basis vector $\bu$ of $K^L$ such that $h(\bu) = 1$ and $\hyp_{xu}$ is a hyperbolic plane, with height given by Lemma \ref{lem_4.7} as
\begin{eqnarray*}
 \H(\hyp_{xu}) & \leq & H(\bx) H(\bu)\\
 & \leq & q^{\frac{(2L^2 - 3L + 2)g(K)}{d}} \left(C_K(L) \E_K(L) H(F)^{\frac{1}{2}} \H(V)\right)^{L-1}.
\end{eqnarray*}

Then, we construct an isotropic vector $\bwy \in K^L$ as in \eqref{y} so that  $\bx,\bwy$ is a hyperbolic pair, and it follows from \eqref{y_height_1} that
\begin{eqnarray}\label{fun2}
h(\bwy) & \leq & H(F_A) h(\bx) h(\bu)^2\nonumber \\
    & =  & H(F_A) h(\bx) \nonumber\\
    & \leq & q^{\frac{(2L^2-3L+2)g(K)}{d}} \left( C_K(L)\E_K(L) H(F)^{\frac{1}{2}} \H(V) \right)^{L+1}.
\end{eqnarray}

Let $U$ be the orthogonal complement of $\hyp_{xu}$ in $(K^L, F_A)$.  Lemma \ref{lem_2.3} implies that
$$\H(U) \leq H(F_A)^2 H(\bx) H(\bu) \leq H(F_A)^2 H(\bx).$$
At the same time, the function field case Siegel's lemma guarantees a small-height basis $\bv_1,\dots,\bv_{L-2}$ for $U$ with
\begin{eqnarray}\label{fun3}
\prod_{i=1}^{L-2} h(\bv_i) & \leq & C_K(L-2) \E_K(L-2)\, \H(U)\nonumber \\
    & \leq & C_K(L) \E_K(L)\, \H(U) \nonumber \\
    & \leq & C_K(L) \E_K(L) H(F_A)^2 H(\bx) \nonumber \\
    & \leq & q^{\frac{(2L^2-3L+2)g(K)}{d}} \left( C_K(L)\E_K(L) \right)^{L+4} \left( H(F)^{\frac{1}{2}} \H(V) \right)^{L+3},
\end{eqnarray}
where the last inequality is from (\ref{fun1}).

Now define the polynomials $Q(\bY)$ and $G(\bY)$ as in the number field case, using the same matrix $B = \left( \bx\ \bwy\ \bv_1 \dots \bv_{L-2} \right)$.  Once again, Lemma~\ref{quad1} guarantees the existence of a point $\bw \in K^L$ such that $Q(\bw) = 0$, $G(\bw) \neq 0$, and
\begin{equation} \label{fun4}
h(\bw) \leq A_K(D + 2)^3 A_K(2 D)^2 H(Q).
\end{equation}
The same procedure of height estimation as in the number field case produces
\begin{eqnarray*}
H(Q) & \leq & H(B^t F_A B) \leq  H(B)^2 H(F_A) \leq H(F_A) h(\bx)^2 h(\bwy)^2 \prod_{i=1}^{L-2} h(\bv_i)^2 \nonumber \\
& \leq & C_K(L)^2 \E_K(L)^{2} H(F) \H(V)^2 h(\bx)^2 h(\bwy)^2 \prod_{i=1}^{L-2} h(\bv_i)^2\\
& \leq & q^{\frac{(12L^2-18L+12)g(K)}{d}} C_K(L)^{6L+10} \E_K(L)^{6L+10} \left( H(F) \H(V)^2 \right)^{3L+4}.
\end{eqnarray*}
As in the number field case, the point $\bz: = A(B\bw)$ is the nontrivial zero of $F$ we desire.  To estimate its height, we first have
\begin{equation*}
h(B\bw) = h \left( w_1 \bx + w_2 \bwy + \sum_{i=1}^{L-2} w_{i+2} \bv_i \right) \leq h(\bw) h(\bx) h(\bwy) \prod_{i=1}^L h(\bv_i).
\end{equation*}
Then combining \eqref{fun1}, \eqref{fun2}, \eqref{fun3}, and \eqref{fun4} leads to \eqref{z_bound_miss}.
\endproof

\begin{rem} \label{poly_rest} Notice that it is also easy to obtain a version of Lemma~\ref{quad1} with a restriction to a subspace $V$ of $K^N$ instead of the whole $K^N$ by applying Lemma \ref{restrict} in the same way as we do it in the proof of Proposition~\ref{miss_hyper_one}.
\end{rem}

Henceforth, we fix a collection $\bs = \{\s_1, \ldots, \s_J\}$ of finite sets of homogeneous polynomials in $K[\bX]$.  For simplicity, we use $M$ to denote the integer $M_{\bs}$ defined for $\bs$ as in \eqref{def_M}.

\begin{cor} \label{cortoprop}
Let the notation be as in the statement of Proposition \ref{miss_hyper_one}.  Suppose that $F$ has a nontrivial zero in $V\setminus\Z_{\bs}$.  Then there must be such a zero $\bx$ such that
\begin{equation}\label{z_bound_miss1}
H(\bx) \leq h(\bx) \leq  T_K(L,M)\, H(F)^{\frac{9L+11}{2}}\, \H(V)^{9L+12}.
\end{equation}
\end{cor}

\proof
For our convenience, let $Z(V,F)$ be the set of nontrivial zeros of $F$ in $V$.  Since $Z(V,F) \nsubseteq \Z_{\bs}$, $Z(V,F) \nsubseteq Z_K(\s_i)$ for all $1 \leq i \leq J$, and so for each $i$ at least one polynomial $P_i$ in $\s_i$ is not identically zero on $Z(V,F)$.  Clearly for each $1 \leq i \leq J$, $Z_K(\s_i) \subseteq Z_K(P_i)$. Define
$$P(\bX) = \prod_{i=1}^J P_{i}(\bX),$$
so that $Z(V,F) \nsubseteq Z_K(P)$ while $\Z_{\bs} \subseteq Z_K(P)$. Then it is sufficient to construct a point of bounded height $\bx \in Z(V,F) \setminus Z_K(P)$.  We may now apply Proposition \ref{miss_hyper_one}.   Notice that since $\deg(P) \leq M$, we have $T_K(L,\deg P) \leq T_K(L,M)$.
\endproof

We are now ready for the proof of Theorem \ref{miss_hyper}.  Recall that for the quadratic space $(V, F)$, $m$ is the dimension of a maximal totally isotropic subspace, which is the sum of the Witt index and the dimension of the radical of $V$.

\proof[Proof of Theorem \ref{miss_hyper}]
Let the notation be as in the statement of Theorem \ref{miss_hyper}. We want to prove the existence of a linearly independent collection of vectors $\bx_1,\dots,\bx_m \in Z(V,F) \setminus \Z_{\bs}$ satisfying \eqref{h_order_1} and \eqref{miss_hyper_bnd}.

Corollary \ref{cortoprop} guarantees the existence of a point $\bx_1 \in Z(V,F) \setminus \Z_{\bs}$ satisfying \eqref{z_bound_miss}. In fact, let $\bx_1$ be a point of smallest possible height in $Z(V,F) \setminus \Z_{\bs}$ satisfying
$$H(\bx_1) \leq h(\bx_1) \leq  T_K(L,M)\, H(F)^{\frac{9L+11}{2}}\, \H(V)^{9L+12}.$$
If $m=1$, we are done; hence suppose that $m > 1$. Then there must exist a maximal totally isotropic subspace $W_1$ of $(V,F)$ containing $\bx_1$, and so $W_1 \nsubseteq \Z_{\bs}$ and $\dim_K W_1 = m$. Then, as implied by Theorem \ref{general_null}, $W_1$ has a full basis $\bu_1,\dots,\bu_m$ outside of $\Z_{\bs}$. Let $\XX_1$ be an $(N-1)$-dimensional subspace of $K^N$ containing $\bx_1$ so that $W_1 \nsubseteq \XX_1$.  Then, at least one of $\bu_1,\dots,\bu_m$ is not in $\XX_1$.  We regard $\XX_1$ as the zero set of a linear form $\bl\in K[\bX]$, and define $\bs^1$ to be the collection containing all the $\s_i$ and $\{\bl\}$; hence $M_{\bs^1} = M + 1$.   Since $W_1 \subseteq Z(V,F)$ but $W_1 \not \subseteq \Z_{\bs^1}$, we can conclude that $Z(V,F) \nsubseteq \Z_{\bs^1}$.  So, we may apply Corollary~\ref{cortoprop} for $\Z_{\bs^1}$ to obtain an $\bx_2 \in Z(V,F) \setminus \Z_{\bs^1}$ such that $\bx_1$ and $\bx_2$ are linearly independent and
$$H(\bx_2) \leq h(\bx_2) \leq T_K(L,M+1)\, H(F)^{\frac{9L+11}{2}}\, \H(V)^{9L+12}.$$
We can assume that $\bx_2$ is a point of smallest possible height in $Z(V,F) \setminus \Z_{\bs^1}$ satisfying the above inequalities. If $m=2$, we are done; otherwise let us assume $m > 2$.  Then there must exist a maximal totally isotropic subspace $W_2$ of $(V,F)$ containing $\bx_1,\bx_2$, and so $W_2 \nsubseteq \Z_{\bs}$ and $\dim_K W_2 = m$.  Again,  Theorem \ref{general_null} guarantees that $W_2$ has a full basis $\bu'_1,\dots,\bu'_m$ outside of $\Z_{\bs}$.  Let $\XX_2$ be an $(N-1)$-dimensional subspace of $V$ containing vectors $\bx_1,\bx_2$ but not the subspace $W_2$, and let $\Z_{\bs^2}$ be the collection containing all the $\s_i$ and the singleton set of a linear form defining $\XX_2$.   Continuing to apply Corollary~\ref{cortoprop}  in the same manner as above, we construct a collection of linearly independent vectors $\bx_1,\dots,\bx_m \in V \setminus \Z_{\bs}$ satisfying
\begin{equation*}
H(\bx_1) \leq H(\bx_2) \leq \dots \leq H(\bx_m),\ h(\bx_1) \leq h(\bx_2) \leq \dots \leq h(\bx_m),
\end{equation*}
and
\begin{equation}
\label{4more}
H(\bx_n) \leq h(\bx_n) \leq T_K(L,M+1)\, H(F)^{\frac{9L+11}{2}}\, \H(V)^{9L+12}.
\end{equation}
This completes the proof of the theorem.
\endproof

We now turn to Corollary~\ref{miss_subspace}. First, we define a constant $a_K(L,N,m)$ by
$$a_K(L,N,m) = \begin{cases}
2^{(2m+1)(L-m-1)} B_K(L-m-1)^{2(L-m-1)} N^2 & \mbox{ for number fields},\\
 q^{\frac{(L-m-1)^2g(K)}{d}} & \mbox{ for function fields}.
\end{cases}$$
Then, we define $T^1_K(L,M,N,m)$, which will appear as the implied constant in \eqref{cor_bnd}, by
\begin{equation*}
\label{cor_bnd_const_nf}
T^1_K(L,M,N,m) = a_K(L,N,m)\, T_K(L,M+1)^2.
\end{equation*}

\proof[Proof of Corollary~\ref{miss_subspace}]
We first show that for each $1 \leq n \leq m$, there exists a maximal totally isotropic subspace $W^m_n$ of $(V,F)$ of bounded height, containing the corresponding point $\bx_n$ from the statement of Theorem~\ref{miss_hyper}; since $\bx_n \notin \Z_{\bs}$, it follows that $W^m_n \nsubseteq \Z_{\bs}$.

First suppose that $\bx_n \in V^{\perp}$, then $\Z_{\bs}$ cannot contain any maximal totally isotropic subspace of $(V,F)$, since each one of them contains $V^{\perp}$. Hence we can pick $W^m_n$ to be a maximal totally isotropic subspace of $(V,F)$ of bounded height as guaranteed by Theorem~1 of \cite{vaaler:smallzeros} in case $K$ is a number field, or by Theorem~\ref{smallspace} above in case $K$ is a function field.  Next assume that $\bx_n$ is a nonsingular zero, then define
$$U_n = \{ \bz \in V : F(\bz,\bx_n) = 0 \} = \{ \bz \in K^N : \bz^t (F \bx_n) = 0 \} \cap V,$$
so that $\dim_K U_n = L-1$ and
\begin{equation}
\label{ht_Un}
\H(U_n) \leq \H(F \bx_n) \H(V) \leq N^{3\delta/2} H(F) H(\bx_n) \H(V),
\end{equation}
by the Brill-Gordan duality principle (discussed in Section~\ref{notation} above), combined with Lemmas~\ref{lem_2.3} and~\ref{intersection} above. Let $W'_n$ be a maximal totally isotropic subspace of $(U_n,F)$ of bounded height as guaranteed by Theorem~1 of \cite{vaaler:smallzeros} in case $K$ is a number field, or by Theorem~\ref{smallspace} above in case $K$ is a global function field.  Therefore
\begin{equation}
\label{ht_max_isot_nf}
\H(W'_n) \leq \begin{cases}
\left( 2^{2m+1} B_K(L-m-1)^2 H(F) \right)^{L-m-1} \H(U_n) & \mbox{ if $K$ is a number field},\\
q^{\frac{(L-m-1)^2g(K)}{d}} H(F)^{L-m-1} \H(U_n) & \mbox{ if $K$ is a function field}.
\end{cases}
\end{equation}

Now, define $W^m_n=\spn_K \{ \bx_n,W'_n \}$, then $W^m_n$ is a maximal totally isotropic subspace of $(V,F)$ containing $\bx_n$. Moreover,
\begin{equation}
\label{5more}
\H(W^m_n) \leq N^{\delta/2}\ H(\bx_n) \H(W'_n).
\end{equation}
Combining this observation with \eqref{ht_max_isot_nf} and the bounds of Theorem~\ref{miss_hyper}, we obtain
$$\H(W^m_n) \leq T^1_K(L,M,N,m)\, H(F)^{10L-m+11} \H(V)^{18L+25},$$
which proves \eqref{cor_bnd}.

Now Siegel's lemma implies the existence of a basis $\bw_1,\dots,\bw_m$ for $W^m_n$ such that
$$\prod_{i=1}^m h(\bw_i) \leq C_K(m) \E_K(m)^{1-\delta} \H(W^m_n).$$
Since $\bo \neq \bx_n \in W^m_n$, there must be a subcollection of $m-1$ of these vectors which are linearly independent with $\bx_n$; since we did not order these vectors by height, we can assume without loss of generality that $\bx_n,\bw_2,\dots,\bw_m$ are linearly independent. Then for each $1 \leq k < m$, define
$$W^k_n = \spn_K \{ \bx_n,\bw_2,\dots,\bw_k \},$$
so that $\bx_n \in W^k_n$, $\dim_K W^k_n = k$,
$$\spn_K \{ \bx_n \} = W^1_n \subset W^2_n \subset \dots \subset W^m_n,$$
By Lemma \ref{lem_4.7}, together with \eqref{4more} and \eqref{5more}, we obtain
\begin{eqnarray*}
\H(W^k_n) & \leq & N^{\delta k/2} H(\bx_n) \prod_{i=2}^k h(\bw_i) \\
    & \leq & N^{\delta k/2} C_K(m) \E_K(m)^{1-\delta} H(\bx_n) \H(W^m_n),
\end{eqnarray*}
and
$$H(\bx_n) \H(W^m_n) \leq  T_K(L,M+1) T^1_K(L,M,N,m) H(F)^{\frac{29L + 33 - 2m}{2}} \H(V)^{27L + 37}.$$
This proves Corollary \ref{miss_subspace}.
\endproof

\appendix
\section{Linear bases of small height}
\label{null_basis}

In this appendix we present two different variations of Siegel's lemma.  Unlike in the rest of the paper, $K$ could be $\qbar$ or a function field with any perfect field for their coefficient fields, i.e., not necessarily just finite fields.  Heights and field constants in this more general situation are defined in the completely analogous manner (see~\cite{null}).

In below, $\bs$ is the same collection of finite sets of homogeneous polynomials in $K[\bX]$ we have been working with, and $M$ is the integer $M_{\bs}$ defined in \eqref{def_M}.

\begin{thm} \label{general_null} Let $K$ be a number field, $\qbar$, or a function field with perfect constant field. Let $N \geq 2$ be an integer, and let $V$ be an $L$-dimensional subspace of $K^N$, $1 \leq L \leq N$. Suppose that $V \nsubseteq \Z_{\bs}$. Then there exists a basis $\bx_1,\dots,\bx_L \in V \setminus \Z_{\bs}$ for $V$ over $K$ such that
\begin{equation}
\label{h_order}
H(\bx_1) \leq H(\bx_2) \leq \dots \leq H(\bx_L),\ h(\bx_1) \leq h(\bx_2) \leq \dots \leq h(\bx_L),
\end{equation}
and for each $1 \leq n \leq L$,
\begin{equation}
\label{gen_bnd_1}
H(\bx_n) \leq h(\bx_n) \leq  L^{\delta} \E_K(L)^{1-\delta} A_K(M+1) C_K(L) \H(V),
\end{equation}
where $\delta$ is as in \eqref{delta}, $C_K(L)$ is as in \eqref{CKL}, $A_K(M+1)$ is as in \eqref{AKLM}, and $\E_K(L)$ is as in \eqref{EKL}.
\end{thm}

\proof
We want to prove the existence of a linearly independent collection of vectors $\bx_1,\dots,\bx_L \in V \setminus \Z_{\bs}$ satisfying \eqref{h_order} and \eqref{gen_bnd_1}. Theorem 1.4 of \cite{null} guarantees the existence of a point $\bx_1 \in V \setminus \Z_{\bs}$ with
\begin{equation}
\label{v1}
H(\bx_1) \leq h(\bx_1) \leq  L^{\delta} \E_K(L)^{1-\delta} A_K(M) C_K(L) \H(V).
\end{equation}
In fact, let $\bx_1$ be a point of smallest possible height in $V \setminus \Z_{\bs}$ satisfying \eqref{v1}. If $L=1$, we are done. If not, let $\XX_1$ be an $(N-1)$-dimensional subspace of $K^N$ containing $\bx_1$ and not containing the entire $V$, and let $\bs^1$ be the collection containing all the $\s_i$ and the singleton set of a linear form defining $\XX_1$.  Then $V \nsubseteq \Z_{\bs^1}$ and $M_{\bs^1} = M + 1$, thus Theorem 1.4 of \cite{null} guarantees the existence of a point $\bx_2 \in V \setminus \Z_{\bs^1}$ with
\begin{equation}
\label{v2}
H(\bx_2) \leq h(\bx_2) \leq  L^{\delta} \E_K(L)^{1-\delta} A_K(M+1) C_K(L) \H(V),
\end{equation}
and we can assume that $\bx_2$ is a point of smallest possible height in $V \setminus \Z_{\bs^1}$ satisfying \eqref{v2}. If $L=2$, we are done. If not, we can let $\XX_2$ be an $(N-1)$-dimensional subspace of $K^N$ containing vectors $\bx_1,\bx_2$  and not containing the entire $V$, and let
$\bs^2$ be the collection containing all the $\s_i$ and the singleton set of a linear form defining $\XX_2$.  Then $V \nsubseteq \Z_{\bs^2}$, and by continuing to apply Theorem 1.4 of \cite{null} in the same manner, we construct a collection of linearly independent vectors $\bx_1,\dots,\bx_L \in V \setminus \Z_{\bs}$ satisfying \eqref{h_order} and \eqref{gen_bnd_1}. This completes the proof of the theorem.
\endproof

\begin{thm} \label{orth_siegel} Let $K$ be a function field with a perfect constant field and $V \subset K^N$ an $L$-dimensional space of $K^N$, $1 \leq L < N$. Let $F$ be a nonzero quadratic form in $N$ variable over $K$, where we also write $F$ for the associated symmetric bilinear form. Then there exists a basis $\bx_1,\dots,\bx_L$ for $V$ over $K$ such that $F(\bx_i,\bx_j) = 0$ for all $i \neq j$, and
$$\prod_{i=1}^L H(\bx_i) \leq C_K(L)^{\frac{L^2+L-2}{4}} H(F)^{\frac{L(L+1)}{2}} H(V)^L,$$
where $C_K(L)$ is as in \eqref{CKL}.
\end{thm}

\noindent
The proof of Theorem~\ref{orth_siegel} is identical to the proof of Theorem 2.4 of \cite{witt}, where we use Thunder's function field version of Siegel's lemma instead of the Bombieri-Vaaler version over a number field.

\bibliographystyle{plain}  
\bibliography{quad_zero-1.6}        

\end{document}